\documentclass[opre,nonblindrev]{informs3aa}

\OneAndAHalfSpacedXI
\usepackage{natbib}
 \bibpunct[, ]{(}{)}{,}{a}{}{,}%
 \def\bibfont{\small}%
 \def\newblock{\ }%

\usepackage{amsmath}
\usepackage{amssymb}
\usepackage{amsfonts}
\usepackage{amsmath}
\usepackage{bbm}
\usepackage[inline]{enumitem}
\usepackage{float}
\usepackage[margin=1in]{geometry}
\usepackage{graphicx}
\usepackage{latexsym}
\usepackage{mathabx}
\usepackage{mathtools}
\usepackage{multirow}
\usepackage{ifpdf}
\ifpdf
    \usepackage[pdftex]{hyperref}
\else
    \usepackage[hypertex]{hyperref}
\fi

\usepackage{xcolor}

\makeatletter
\newcommand{\eqnum}{\refstepcounter{equation}\textup{\tagform@{\theequation}}}
\makeatother

\newcommand{\R}{\mathbb{R}}

\newcommand{\E}[1]{{\mathbb E}\left[#1\right]}

\newcommand{\1}[1]{{\mathbbm 1}\left\{ #1 \right\}}
\newcommand{\F}[1]{{F}(#1)}

\newcommand\cG{\mathcal G}

\newcommand\cI{\mathcal I}

\newcommand\cW{\mathcal W}

\newcommand{\ceil}[1]{\left\lceil\, #1 \,\right\rceil}

\let\oldmarginpar\marginpar
\renewcommand\marginpar[1]{\-\oldmarginpar[\raggedleft\scriptsize #1]%
{\raggedright\scriptsize #1}}

\TheoremsNumberedThrough     % Preferred (Theorem 1, Lemma 1, Theorem 2)
\ECRepeatTheorems

\EquationsNumberedThrough    % Default: (1), (2), ...
\begin{document}

\TITLE{{Logarithmic Regret in the \\Dynamic and Stochastic Knapsack Problem \\with Equal Rewards}}

\RUNTITLE{Logarithmic regret in the dynamic and stochastic knapsack problem with equal rewards}

\ARTICLEAUTHORS{%
\AUTHOR{Alessandro Arlotto}
\AFF{The Fuqua School of Business, Duke University, 100 Fuqua Drive, Durham, NC, 27708, \\  \EMAIL{alessandro.arlotto@duke.edu}}

%\AUTHOR{Yehua Wei}
%\AFF{Carroll School of Management, Boston College, 140 Commonwealth Avenue, Chestnut Hill, MA 02467,\\ \EMAIL{yehua.wei@bc.edu}}

\AUTHOR{Xinchang Xie}
\AFF{The Fuqua School of Business, Duke University, 100 Fuqua Drive, Durham, NC, 27708, \\  \EMAIL{xinchang.xie@duke.edu}}
}

\RUNAUTHOR{A. Arlotto and X. Xie}

\ABSTRACT{
We study a dynamic and stochastic knapsack problem in which a decision maker is sequentially presented
with items arriving according to a Bernoulli process
over $n$ discrete time periods.
Items have equal rewards and independent weights that are drawn from a known non-negative continuous distribution $F$.
The decision maker seeks to maximize the expected total reward of the items that she includes in the knapsack
while satisfying a capacity constraint and while  making terminal decisions as soon as each item weight is revealed.
Under mild regularity conditions on the weight distribution $F$, we prove that the regret---the expected
difference between the performance of the best sequential algorithm and that of a prophet who
sees all of the weights before making any decision---is, at most, logarithmic in $n$.
Our proof is constructive. We devise a reoptimized heuristic that achieves this regret bound.
}


\MSCCLASS{Primary:    90C39. %dynamic programming
    Secondary:  60C05, %Combinatorial probability
                68W27, %Online algorithms
                68W40, %Analysis of algorithms
                90C27.
}

\ORMSCLASS{Primary:
        Analysis of Algorithms: Suboptimal algorithms.
        Dynamic programming: Markov.
    Secondary:
        Decision Analysis: Sequential, Theory.
        %Probability: Markov processes;
}

\HISTORY{
		First version: September 6, 2018. \emph{This version:} October 28, 2019.%\today. %June 20, 2019. 	\\
		%\emph{File name:} \texttt{\jobname.tex}
}

\maketitle

\section{Introduction}\label{se:introduction}

The knapsack problem is one of the classic problems in operations research.
It arises in resource allocation, and it counts
numerous applications in auctions, logistics, portfolio optimization, scheduling, and transportation among others
\citep[cf.][]{MarTot:WILEY1990,KellererPferschyPisinger:SPRINGER2004}.
In its dynamic and stochastic formulation \citep[see, e.g.][]{PapRajKle:MS1996,KlePap:OR1998,KlePap:OR2001}
a decision maker (referred to as \emph{she}) is given a knapsack with finite capacity $0 \leq c < \infty$
and is sequentially presented with items arriving over a time horizon with $n$ discrete time periods,
indexed by $i \in [n] \equiv \{1,2, \ldots, n\}$.
In each period $i\in[n]$, an item arrives with probability $p$,
its weight-reward pair $(\mathfrak W_i,\mathfrak R_i)$ is revealed,
and the decision maker needs to decide whether to include the arriving item in the knapsack
or to reject it forever.
Here, the weight $\mathfrak W_i$ represents the amount of knapsack capacity that the item arriving in period $i$
consumes if the decision maker chooses to include it in the knapsack,
and the reward $\mathfrak R_i$ represents what the decision maker collects upon inclusion.
The pairs $(\mathfrak W_i, \mathfrak R_i)$, $i \in [n]$, are independent and with common, known, bivariate distribution supported on the
nonnegative orthant.

By imposing different assumptions on the weight-reward distribution,
one recovers knapsack instances of independent interest.
For instance, in the problem of real-time uniprocessor scheduling under conditions of overload \citep[see, e.g.,][]{BaruahHaritsaSharma:1994}
a decision maker wants to maximize the number of jobs that are processed on a single machine by a fixed deadline.
In this context, the deadline is the knapsack capacity and jobs correspond to items.
Their rewards are all equal to one, and their durations
correspond to the item weights.
This scheduling application motivates the model in this paper.
We assume that the rewards are deterministic and all equal\footnote{This also covers random rewards with common distribution that are revealed only after the inclusion decision.} to $r>0$,
and the weights are independent random variables with common continuous distribution $F$.
We model item arrivals by considering a Bernoulli process $\mathfrak B_1,\mathfrak B_2,\ldots,\mathfrak B_n$
that is independent of everything else, and that is given by a sequence of independent Bernoulli random variables
with success probability $p$.
We then equivalently redefine the weight distribution so that
a no arrival corresponds to the arrival of an item with arbitrarily large weight.
That is, we assume that an item arrives in each period $i\in[n]$
and that its weight is given by the random variable $W_i$ defined by
$$
W_i = \begin{cases}
        +\infty          & \mbox{if } \mathfrak B_i = 0 \\
        \mathfrak{W}_i   & \mbox{if } \mathfrak B_i = 1.
      \end{cases}
$$
We say that a policy $\pi$ is \emph{feasible} if the sum of the weights of the items selected by $\pi$  does not
exceed the knapsack capacity $c$, and we say that the policy is \emph{online} (or \emph{sequential})
if the decision to select item $i$ with weight $W_i$ depends only on the information available up to
and including time $i$.
We then let $\Pi(n,c,p)$ be the set of feasible online policies, and we compare the performance
of the best online policy to that of a prophet who has full (or \emph{offline}) knowledge of the weights $W_1, W_2, \ldots, W_n$
before making any selection.
Under some mild technical conditions on the weight distribution $F$,
we prove that the regret---the expected gap between the performance of the best online policy and its offline counterpart---is
bounded by the logarithm of $n$.
Our proof is constructive. We propose a reoptimized heuristic that exhibits logarithmic regret.
The heuristic is based on resolving  some related optimization problem
at any given time $i \in [n]$ by using the current---rather than the initial---level of remaining capacity as constraint.
The solution of this optimization problem provides us with a state- and time-dependent
threshold that mimics that of the optimal online policy.

If all of the weights $W_1, W_2, \ldots, W_n$ are revealed to the decision maker before she makes any selection,
then her choice is obvious.
To maximize the total reward she collects,
she just sorts the items according to their weights and
selects them starting from the smallest weight and continuing until the knapsack capacity is exhausted.
Formally, if $W_{(1,n)} \leq W_{(2,n)} \leq \cdots \leq W_{(n,n)}$ are the order statistics of $W_1, W_2, \ldots, W_n$,
then the maximal reward $R^*_n(c,p,r)$ that the decision maker collects is given by
\begin{equation}\label{eq:offline-sort-algorithm}
R^*_n(c,p,r)
= \max \left\{ rm: ~m \in \{0,1,\ldots,n\} \mbox{ and }\sum_{\ell = 1}^{m} W_{(\ell, n)} \leq c \right\}.
\end{equation}

Here we compare the total reward of the offline-sort algorithm  \eqref{eq:offline-sort-algorithm}, $R^*_n(c,p,r)$, with that
of an online feasible policy $\widehat \pi \in \Pi(n,c,p)$
that is based on a sequence of reoptimized time- and state-dependent threshold functions
$\widehat{h}_n,\widehat{h}_{n-1},\ldots,\widehat{h}_1$.
If the current level of remaining capacity is $x$
and the weight of item $i$ is about to be revealed,
then the decision maker computes the threshold $\widehat h_{n-i+1}:[0,\infty) \rightarrow [0,\infty)$
such that $\widehat h_{n-i+1}(x) \leq x$, and she selects item $i$ if and only if
the weight $W_i \leq \widehat h_{n-i+1}(x)$.
Thus if $\widehat{X}_0 = c$ and for $i \in [n]$ one defines the remaining
capacity process $\widehat{X}_i$ recursively by
\begin{equation*}%\label{eq:remaining-capacity-MC}
\widehat{X}_i = \begin{cases}
    \widehat{X}_{i-1}             & \text{if } W_i > \widehat{h}_{n-i+1}(\widehat{X}_{i-1}) \\
    \widehat{X}_{i-1} - W_i       & \text{if } W_i \leq \widehat{h}_{n-i+1}(\widehat{X}_{i-1}),
\end{cases}
\end{equation*}
then the total reward collected by the reoptimized policy $\widehat{\pi}$ can be written as
$$
  R_n^{\widehat \pi} (c,p,r) = \sum_{i=1}^n r \1{W_i \leq \widehat{h}_{n-i+1}(\widehat{X}_{i-1})}.
$$

The random variables $R^*_n(c,p,r)$ and $R_n^{\widehat{\pi}} (c,p,r)$ crucially depend on the weight distribution $F$.
This dependence is mostly expressed through a \emph{consumption function}
$\epsilon_{k p}: [0, \infty) \rightarrow [0,\infty]$ that is defined for $p\in(0,1]$ and for all $1\leq k <\infty$ by
\begin{equation}\label{eq:epsilon-k-definition}
  \epsilon_{k p}(x) = \sup\bigg\{ \epsilon \in [0, \infty) : \int_{0}^{\epsilon} w \, d\F{w} \leq \frac{x}{k p}\bigg\}.
\end{equation}
The consumption function depends on two quantities.
The argument $x$ that denotes the current level of remaining capacity of the knapsack,
and the index $k p$ that refers to the expected number of items with $F$-distributed weights (or \emph{arrivals})
that are yet to be presented to the decision maker.
Furthermore, the function $\epsilon_{k p}(x)$ is always well defined.
If $\mu = \E{\mathfrak W_1} = \E{W_1 \mid B_1 = 1}$ and $k p \mu < x < \infty$ then $\epsilon_{k p}(x) = +\infty$.
Otherwise, the value $\epsilon_{k p}(x)$ satisfies
the integral representation
\begin{equation}\label{eq:epsilon-k-integral}
  \int_0^{\epsilon_{k p}(x)} w   \, d\F{w} = \frac{x}{k p}
  \quad \quad \text{for all }  x \in [0, k p\mu].
\end{equation}
The representation \eqref{eq:epsilon-k-integral} offers an important insight regarding the role
of the consumption function $\epsilon_{k p}(x)$.
The integral on the left-hand side is the \emph{expected reduction}
in the remaining capacity of the knapsack when the current level of remaining capacity is equal to $x$,
and the decision maker selects an
item with weight smaller than $\epsilon_{k p}(x)$.
The function $\epsilon_{k p}(x)$ is then defined  so that
the expected reduction in capacity is equal to the ratio of the current capacity, $x$,
to the expected number of remaining arrivals, $kp$.
That is, the threshold $\epsilon_{k p}(x)$ is constructed so that---in expectation---the available capacity
is spread equally over the remaining arrivals.

As we will see shortly, the threshold $\epsilon_{k p}(x)$  drives most of the estimates in this paper
and, together with the continuity of the weight distribution $F$, it immediately provides us with
an easy upper bound for $\E{R^*_n(c,p,r)}$.
The same threshold together with some mild regularity conditions on the weight distribution $F$
also drives the lower bound for $\E{R_n^{\widehat{\pi}}(c,p,r)}$.
The class of weight distributions we consider for the lower bound is characterized in the next definition.

\begin{definition}[Typical class of distributions with continuous density]\label{def:typical-class}
We say that a non-negative distribution $F$ with continuous density function $f$
belongs to the \emph{typical class} if for some $\bar{w} > 0$,
the following two conditions hold.
\begin{enumerate}[label=(\roman*)]%[align=left,labelindent=!,labelwidth=!,labelsep=0pt]
    \item \textsc{Behavior at zero.}
        There are $0 < \lambda < 1$ and $0 < \gamma < 1$ such that
        \begin{equation}\label{eq:cdf-condition}
        \frac{\F{\lambda w}}{\F{w}} \leq \gamma < 1  \quad \quad \text{for all $w \in (0, \bar{w})$.}
        \end{equation}

    \item \textsc{Monotonicity.}
        The map $w \mapsto w^3 f(w)$ is non-decreasing on $(0, \bar{w})$. That is,
        \begin{equation}\label{eq:monotonicity-wf(w)}
            w_1^3 f(w_1) \leq w_2^3 f(w_2)
            \quad \quad \text{for all $0 < w_1 \leq w_2 < \bar{w}$.}
        \end{equation}
\end{enumerate}
\end{definition}

The class of typical distributions is wide enough to include most well-known non-negative distributions.
In Section \ref{se:typical-class}, we provide specific examples as well as class properties,
but for now we emphasize that the breadth of the typical class comes from the role of the
distribution-dependent parameter $\bar{w}>0$.
Conditions \eqref{eq:cdf-condition} and \eqref{eq:monotonicity-wf(w)}
need only to hold near zero---or, more precisely, on $(0,\bar w)$---and not on the full support of the weight distribution or on the whole capacity interval $[0,c]$.
In fact, for many distributions the parameter $\bar w$ for which
\eqref{eq:cdf-condition} and \eqref{eq:monotonicity-wf(w)} hold is much smaller than
the minimum between the initial capacity and the supremum of the support.

The main results of this paper are gathered in the theorem below.
First, we provide an upper bound for $\E{R^*_n(c,p,r)}$ that holds
for any continuous distribution $F$.
Then, we turn to distributions that belong to the typical class,
and we prove that there is a matching lower bound.
As a by-product of our analysis, we establish that
the regret is, at most, $O(\log n)$ as $n \rightarrow \infty$.\footnote{Throughout this paper, the function $\log$ denotes the natural logarithm.}
While our theoretical result provides only a regret bound,
related results and the numerical experiments of Section \ref{se:numerical}
tell us  that the regret bound is actually of the correct order.

\begin{theorem}[Logarithmic regret bound]\label{thm:main}
Consider a knapsack problem with capacity $0 \leq c < \infty$ and
with items that arrive over $1 \leq n <\infty$ periods
according to a Bernoulli process with arrival probability $p\in(0,1]$.
If the items have rewards equal to $r$ and weights with continuous distribution $F$,
then
$$
\max_{\pi \in \Pi(n,c,p)} \E{ R_n^{\pi} (c,p,r) } \leq \E{ R^*_n(c,p,r) } \leq  npr \F{ \epsilon_{n p}(c) }.
$$
Furthermore, there is a feasible online policy $\widehat \pi \in \Pi(n,c,p)$ such that
if the weights are independent and their distribution $F$ belongs to the typical class then there is a constant $1 < M < \infty$
depending only on $F$, $p$, and $r$ for which
$$
 npr \F{ \epsilon_{n p}(c) } - M( 1+  \log n ) \leq \E{ R_n^{\widehat{\pi}} (c,p,r) } \leq \max_{\pi \in \Pi(n,c,p)} \E{ R_n^{\pi} (c,p,r) }.
$$
In turn, if the weights are independent and the distribution $F$ belongs to the typical class,
then we have the regret bound
$$
\E{ R^*_n(c,p,r) } - \max_{\pi \in \Pi(n,c,p)} \E{ R^\pi_n(c,p,r)  } \leq \E{ R^*_n(c,p,r) } - \E{ R^{\widehat \pi}_n (c,p,r) } \leq M ( 1 + \log n ).
$$
\end{theorem}

The special case with deterministic arrivals and unitary rewards has been extensively studied in the literature.
The upper bound $\E{ R^*_n(c,1,1) } \leq  n \F{ \epsilon_n(c) }$ was first proved by \citet{BruRob:AAP1991}.
Here, we provide a generalization that is based on a relaxation of some appropriate optimization problem.
The solution to this relaxation is the basis for constructing the reoptimized heuristics $\widehat \pi$.
The lower bound $\E{ R_n^{\widehat{\pi}} (c,p,r) } \geq npr \F{ \epsilon_{n p}(c) } - O( \log n )$ as $n \rightarrow \infty$
is essentially new, and it substantially improves on existing estimates.
The best results to date for general weight distribution $F$ are due to \citet{RheTal:JAP1991}
who study a non-adaptive heuristic and prove that
\begin{equation}\label{eq:RheeTalagrand-LowerBound}
n\F{\epsilon_n(c)}\left\{ 1- \left[\frac{\epsilon_n(c)}{c}\right]^{1/2} - \frac{\epsilon_n(c)}{c} \right\}
\leq \max_{\pi \in \Pi(n,c,1)} \E{ R_n^{\pi} (c,1,1) }
\quad \quad \text{for all } n \geq 1.
\end{equation}
For instance, if $F(x) = \sqrt{x}$ for $x\in (0,1)$ then the lower bound \eqref{eq:RheeTalagrand-LowerBound} implies
an upper bound for the regret that is $O(n^{1/3})$ as $n \rightarrow \infty$.
Similarly, if $F(x) = x^2$ for $x \in (0,1)$ then the same lower bound gives us a regret
upper bound that behaves like $O(n^{1/6})$ as $n \rightarrow \infty$.

A case that deserves special attention is when $F$ is the uniform distribution on the unit interval, the reward $r=1$,
and the initial capacity $c=1$.
In this context, the \citet{RheTal:JAP1991} lower bound provides us with a regret upper bound that
behaves like $O(n^{1/4})$ as $n\rightarrow \infty$, but better bounds are available in the literature.
This special dynamic and stochastic knapsack problem is in fact equivalent to the problem of the
sequential selection of a monotone decreasing subsequence from a sample of $n$ independent observation
with the uniform distribution on the unit interval \citep[cf.][]{SamSte:AP1981}.
The equivalence was first observed by \citet[][pp. 457--458]{CofFlaWeb:AAP1987},
and it can be established by observing that the Bellman equations for the two problems are
the same after a change of variable.
Informally, if the number of remaining periods is the same in both problems
and the current capacity of the knapsack is equal to the last selected subsequence element,
then the largest weight that is optimal for inclusion is equal to the maximum amount the decision maker is
willing to go down in optimally selecting a new subsequence element.
Since the weights as well as the subsequence elements are both uniformly distributed on the unit interval,
these two actions happen with the same probability.
For this subsequence-selection problem, \citet{ArlNguSte:SPA2015,ArlottoWeiXie:RSA2018}
prove that the expected performance $\nu^*_n$ of the best online policy satisfies the estimate
$\nu^*_n = \sqrt{2n} - O(\log n)$ as $n \rightarrow \infty$.
The equivalence between the two problems, however, holds \emph{only} for uniform weights.
As Theorem \ref{thm:main} suggests, the weight distribution $F$ plays a crucial role in the
estimates for the dynamic and stochastic knapsack problem with equal rewards.
Instead, the monotone subsequence problem is distribution invariant,
and one can consider uniformly distributed subsequence elements without loss of generality.
More importantly, \citet{Seksenbayev:WP2018} and \citet{GnedinSeksenbayev:WP2019} characterize the second order asymptotic expansion of $\nu^*_n$ and
establish that $\nu^*_n = \sqrt{2n} - \tfrac{1}{12} \log n + O(1)$ as $n\rightarrow \infty$.
This remarkable result tells us that our regret bound is order tight, and that no
online algorithm can---at this level of generality---be within $O(1)$ of offline sort.

\subsection*{Organization of the paper}

The paper is organized as follows.
In Section \ref{se:literature-review}, we review the related literature.
In Section \ref{se:prophet-upper-bound}, we prove the prophet upper bound $\E{R^*_n(c,p,r)} \leq npr \F{\epsilon_{n p}(c)}$
by showing that the offline-sort algorithm  \eqref{eq:offline-sort-algorithm}
can be reinterpreted as a parsimonious threshold policy
and by solving a relaxation of some related optimization problem.
This solution then guides us in the construction of policy $\widehat \pi$
that is presented in Section \ref{se:heuristic-policy}.
In Section \ref{se:typical-class}, we discuss the generality of the typical class of distributions,
and we derive some properties that we then use---in Section \ref{se:main-proof}---to
prove that the reoptimized policy $\widehat \pi$ exhibits logarithmic regret.
In Section \ref{se:numerical}, we present numerical experiments that
provide further insights into our regret bound, while
in Section \ref{se:weights-w-multiple-types} we discuss weight distributions with multiple types.
Finally, in Section \ref{se:conclusion} we make closing remarks and underscore some open problems.

\section{Literature review: knapsack problems and approximations}\label{se:literature-review}

Knapsack problems uniquely combine simple formulations, non-trivial mathematical analyses,
and relevance in several application-driven domains.
As such, different knapsack problems have been considered in the literature,
and a lot of effort has been devoted to the development of (near-) optimal policies.
Most of the differences that have been accounted for concern
the item arrival process (\emph{static} versus \emph{dynamic}),
the probabilistic assumptions on the weight-reward pairs (deterministic and/or stochastic),
and the objective of the decision maker (reward maximization, target achievement, etc.).

For instance, in the early formulation of \citet{Dantzig:OR1957}, we have a \emph{static} model with a finite number of items
that are all available before any decision is made and have deterministic weights and deterministic rewards.
The decision maker then seeks to find a maximum-reward subset of these items with total weight that does not
exceed a capacity constraint.
Following this classic formulation, researchers have considered several \emph{static} knapsack instances with randomness
in the weights and/or in the rewards.
While studying a scheduling problem, \citet{DermanLiebermanRoss:MS1978} studied a \emph{static} and stochastic knapsack problem
with items that belong to different categories.
Items that belong to the same category have common deterministic rewards
and independent, exponentially distributed weights with category-dependent parameter.
The decision maker then seeks to maximize total expected rewards when
the realized weights are revealed only after each item is included in the knapsack.
The authors prove that the greedy policy based on reward-to-mean-weight ratios is optimal.
Analogous \emph{static} and stochastic knapsack problems have been considered by several authors,
including \citet{DeanGoemansVondrak:FOCS2004,DeanGoemansVondrak:SDA2005,DeanGoemansVondrak:MOOR2008},
\citet{BhalgatGoelKhanna:SDA2011}, \citet{LiYuan:STOC2013}, \citet{BladoHuToriello:SIAM2016}, \citet{Ma:MOR2018}, \citet{BladoToriello:SIAM2019},
and \citet{BalseiroBrown:OR2019}.
\citet{GuptaKrishnaswamyMolinaroRavi:FOCS2011} and \citet{MerzifonluogluGeunesRomeijn:MPROG2012}
follow along similar lines, but consider both random weights and random rewards.
Most notably, \citet{DeanGoemansVondrak:FOCS2004,DeanGoemansVondrak:SDA2005,DeanGoemansVondrak:MOOR2008}
study a \emph{static} and stochastic knapsack problem with deterministic rewards and independent random weights with
arbitrary distributions that are realized only upon insertion in the knapsack.
They construct a polynomial time adaptive policy that is within a constant multiplicative gap,
and they compare the performance of adaptive and non-adaptive policies.
Their work is particularly relevant to us as it is among the first ones to assess the benefits of adaptivity.

\emph{Static} stochastic knapsack problems have also been studied under different optimization objectives.
For instance, there is a stream of related literature that considers \emph{static} stochastic knapsack problems
(typically with deterministic weights and random rewards) in which the objective is to
maximize the probability that the total reward will achieve a certain given target.
\citep[See, e.g.,][among others.]{Henig:OR1990,CarrawaySchmidtWeatherford:NRL1993,IlhanIravaniDaskin:OR2011}

Alongside the static knapsack problems mentioned thus far there are several \emph{dynamic} models
in which items arrive over time and their weight-reward pairs
are revealed to the decision maker who irrevocably decides on inclusion in the knapsack as soon as each item arrives and
without seeing the weights and/or the rewards of future items.
\emph{Dynamic} and stochastic knapsack problems are widespread.
For instance, if one assumes that the weights are all equal to one and that
the rewards are random, then one recovers the multi-secretary problem
\citep[see, e.g.][]{Cay:ET1875,Mos:SM1956,kleinberg2005multiple}.
For this problem, \citet{ArlottoGurvich:SSY2019} prove that if the reward distribution is discrete, then the regret is uniformly bounded in
the number of items and the knapsack capacity.
Similarly, if one assumes that the rewards are all equal to one and that the weights are random,
then one finds an instance of the single-machine scheduling problem of \citet{BaruahHaritsaSharma:1994} that motivates this paper.
Finally, when both the weights and the rewards are random, one recovers---among others---the
sequential investment problems of \citet{DermanLiebermanRoss:OPRE1975} and \citet{Prastacos:MANSCI1983},
or the multi-secretary problem of \citet{Nakai:OPRE1986} which allows for an unknown number of applicants in each period.
When both the weights and the rewards are random, few regret bounds are available.
A notable exception is the work of \citet{Marchetti-SpaccamelaVercellis:MP1995} who prove a $O(\log^{3/2} n)$  regret bound
when both the weights and the rewards are independent and uniformly distributed on the unit interval, and
the knapsack capacity is proportional to the number of periods.
For the same formulation, \citet{Lueker:JOFALG1998} improves
\citeauthor{Marchetti-SpaccamelaVercellis:MP1995}'s result to $O(\log n)$ and shows that it is best possible.

Multi-dimensional generalizations of the \emph{dynamic} and stochastic knapsack problem have found several
applications in revenue management and resource allocation.
In the network revenue management problem, heterogeneous customers belonging to different classes
arrive sequentially over time, request a product, and offer a price.
If the request is accepted, then a collection of resources that constitute the product
is depleted, and the offered price is earned.
Otherwise the resource capacities remain unchanged and the offered price is lost \citep[cf.][]{GallegoVanRyzin:OPRE1997,TallvanR:KLU2004}.
The solution of the network revenue management problem is famously difficult,
and scholars have studied several non-adaptive as well as adaptive heuristics and proved regret bounds.
A classic non-adaptive approximation scheme based on a deterministic linear-programming relaxation
was studied by \citet{GallegoVanRyzin:MANSCI1994,GallegoVanRyzin:OPRE1997}.
In contrast, adaptive policies have been considered by allowing for periodic reoptimization.
Despite a few specific negative results by \citet{Cooper:OPRE2002}, \citet{ChenHomenDeMello:AOR2010}, and \citet{JasinKumar:OPRE2013},
there are ways to construct reoptimized policies that perform well.
For instance, \citet{ReimanWang:MOOR2008} propose a probabilistic allocation rule that works well with one
reoptimization instance.
\citet{JasinKumar:OPRE2012} and \citet{WuSrikantLiuJiang:NIPS2015} consider a probabilistic allocation
rule that is based on reoptimizing in every period and show that it exhibits uniformly bounded regret
provided that the optimal solution to the original deterministic linear programming relaxation is non-degenerate.
\citet{BumpensantiWang:MS2019} and \citet{VeraBanerjee:WP2018} prove that the uniform regret
bound holds in general, without the non-degeneracy assumption.

\section{A prophet upper bound}\label{se:prophet-upper-bound}

The performance of any online algorithm is bounded above by the full-information (or offline) sort.
If the decision maker knows all of the weights $W_1, W_2, \ldots, W_n$
before making any decision,
then the total reward she collects is the largest number $rm$
such that the sum of the smallest $m$ realizations
does not exceed the capacity constraint.
That is, if $W_{(1,n)} \leq W_{(2,n)} \leq \cdots \leq W_{(n,n)}$
are the order statistics of $\cW \equiv \{W_1, W_2, \ldots, W_n\}$,
then the total reward $R^*_n(c,p,r)$ of offline selections when the initial knapsack capacity is $c$
and the arrival probability is $p$
is given by
\begin{equation}\label{eq:offline-sort}
  R^*_n(c,p,r) = \max \left\{ rm: ~m\in \{0,1,\ldots,n\}, \sum_{\ell = 1}^{m} W_{(\ell, n)} \leq c \text{ and } W_{(\ell, n)} \in \cW \text{ for all } \ell \in [n] \right\}.
\end{equation}
Earlier work has considered unitary rewards and deterministic arrivals by studying the
random variable $R^*_n(c,1,1)$.
First along this line of research, \citet{CofFlaWeb:AAP1987} showed that
$$
R^*_n(c,1,1)  \sim n \F{ \epsilon_n(c) }
\text{ in probability as }
n \rightarrow \infty,
$$
provided that the weight distribution $F$ is continuous, strictly increasing in $w$ when $\F{w}< 1$,
and $\F{ w } \sim A w^\alpha$ as $w \rightarrow 0$ for some $A, \alpha>0$.
Four years later, \citet{BruRob:AAP1991} proved that the same result holds under more general conditions,
and \citet{boshuizen1999smallest} established the asymptotic normality
of $R^*_n(c,1,1)$ after the usual centering and scaling
for different classes of weight distribution $F$.
Lemma 4.1 in \citet{BruRob:AAP1991} is particularly relevant to our discussion here since it tells us that
$$
\E{ R^*_n(c,1,1) } \leq n \F{ \epsilon_n(c) }
\qquad \qquad \text{ for all } n \geq 1.
$$

Here, we generalize this result by accounting for Bernoulli arrivals with probability $p\in(0,1]$
and rewards equal to $r>0$.
Specifically, we show that
$$
\E{ R^*_n(c,p,r) } \leq npr \F{ \epsilon_{n p}(c) }
\qquad \qquad \text{ for all } n \geq 1.
$$

Our proof relies on the observation that the offline-sort algorithm  \eqref{eq:offline-sort}
can be equivalently described as an algorithm
that selects items with weight that is below some threshold.
For any given realization $W_1, W_2, \ldots, W_n$,
the offline-sort algorithm  selects $N^*_n \equiv R^*_n(c,p,1)$ items so
one can compute the value $W_{ (N^*_n, n)}$ of the largest
weight that is selected for inclusion,
and one can then select all of the items $i \in [n]$ that have weight $W_i \leq W_{ (N^*_n, n)}$.
A shortcoming of this interpretation is that one needs to know the
realization of the weight $W_i$ (as well as the realizations of all of the other weights)
to compute the threshold $W_{ (N^*_n, n)}$.
As it turns out, this is not needed in general.
The next lemma shows that there is a thresholding algorithm
that makes the same selections of offline sort,
but in which the threshold used to decide whether to select an item is computed
without using the information about that item's weight.

\begin{lemma}[Threshold policy equivalence]\label{lm:threshold-policy-equivalence-offline-sort}
Let $W_{(1,n)} \leq W_{(2,n)} \leq \cdots \leq W_{(n,n)}$
be the order statistics of $\cW \equiv \{W_1, W_2, \ldots, W_n\}$
and, for $i \in [n]$, let $W_{(1,n-1)} \leq W_{(2,n-1)} \leq \cdots \leq W_{(n-1,n-1)}$
be the order statistics of $\cW_i = \cW \backslash \{W_i\}$.
Then, for
\begin{equation}\label{eq:tau-i}
\tau^i_{n-1} = \max \left\{ m \in \{0,1,\ldots,n-1\}:~\sum_{\ell = 1}^{m} W_{(\ell, n-1)} \leq c \text{ and } W_{(\ell, n-1)} \in \cW_i \text{ for all } \ell \in [n-1] \right\}
\end{equation}
and $N^*_n \equiv R^*_n(c,p,1)$,
we have that
\begin{equation}\label{eq:threshold-equivalence}
    W_i \leq W_{(N^*_n, n)}
    \qquad \text{if and only if}\qquad
    W_i \leq h(\cW_i) \equiv \max\big\{ W_{(\tau^i_{n-1}, n-1)}, c - \sum_{\ell=1}^{\tau^i_{n-1}} W_{(\ell, n-1)}\big\}.
\end{equation}
In turn, it follows that
\begin{equation}\label{eq:Nstar-threhsold-representation}
  R^*_n(c,p,r) = \sum_{i=1}^{n} r\1{ W_i \leq h(\cW_i) }.
\end{equation}
\end{lemma}

\proof{Proof of Lemma \ref{lm:threshold-policy-equivalence-offline-sort}.}
The equivalence \eqref{eq:Nstar-threhsold-representation} is an obvious consequence of \eqref{eq:threshold-equivalence},
so we focus on proving the latter.
If $N^*_n = n$ we have that $\tau^i_{n-1} = n-1$ and $W_i \leq c - \sum_{\ell=1}^{\tau^i_{n-1}} W_{(\ell, n-1)}$
for all $i \in [n]$, so equivalence \eqref{eq:threshold-equivalence} immediately follows.
Instead, if $N^*_n < n$ the proof of \eqref{eq:threshold-equivalence} requires more work.
As a warm-up we note that since the sets $\cW$ and $\cW_i$ differ only in one element, then
\begin{equation}\label{eq:order-statistic-skip-one}
W_{ (\ell, n) } \leq W_{ (\ell, n-1) } \leq W_{ (\ell+1, n) } \qquad \text{ for all } \ell \in [n-1].
\end{equation}
If we now recall the definitions of $\tau^i_{n-1}$ and $N^*_n$
and use the inequalities above we obtain that
$$
\sum_{\ell=1}^{\tau^i_{n-1}} W_{ (\ell, n) } \leq \sum_{\ell=1}^{\tau^i_{n-1}} W_{ (\ell, n-1) } \leq c
\qquad \text{and} \qquad
\sum_{\ell=1}^{N^*_n-1} W_{ (\ell, n-1) } \leq \sum_{\ell=1}^{N^*_n-1} W_{ (\ell+1, n) }  \leq \sum_{\ell=1}^{N^*_n} W_{ (\ell, n) } \leq c.
$$
These two bounds respectively tell us that the offline-sort algorithm  on $\cW$  selects at least $\tau^i_{n-1}$ observations,
and that the same algorithm on $\cW_i$ selects at least $N^*_n - 1$ items.
Thus, it follows that
$$
N^*_n - 1 \leq \tau^i_{n-1} \leq N^*_n,
$$
and we use these bounds to  prove the equivalence \eqref{eq:threshold-equivalence}.

\emph{If.} We now suppose that $W_i \leq h(\cW_i) \equiv \max\big\{ W_{(\tau^i_{n-1}, n-1)}, c - \sum_{\ell=1}^{\tau^i_{n-1}} W_{(\ell, n-1)}\big\}$,
and we seek to show that $W_i \leq W_{ (N^*_n, n) }$.
We consider two cases, one per each possible realization of $\tau^i_{n-1}$.
\begin{description}
    \item[\textsc{Case 1:} $\tau^i_{n-1} = N^*_n - 1$.]
        If $\tau^i_{n-1} = N^*_n - 1$ then the definition of $\tau^i_{n-1}$
        in \eqref{eq:tau-i} tells us that
        $$
        c - \sum_{\ell=1}^{N^*_n - 1} W_{ (\ell, n-1) } < W_{ (N^*_n, n-1) },
        $$
        so if we apply the right inequality of \eqref{eq:order-statistic-skip-one} to $\ell = N^*_n-1$ and $\ell = N^*_n$, we obtain that
        \begin{equation}\label{eq:maximands-bounds}
        W_{ (N^*_n-1, n-1) } \leq W_{ (N^*_n, n) }
        \qquad\qquad \text{and} \qquad\qquad
        c - \sum_{\ell=1}^{N^*_n - 1} W_{ (\ell, n-1) } < W_{ (N^*_n+1, n) }.
        \end{equation}
        If $W_{ (N^*_n, n) } = W_{ (N^*_n + 1, n) }$ then the two inequalities in \eqref{eq:maximands-bounds}
        give us that
        $
        h(\cW_i) = \max\big\{ W_{(N^*_n - 1, n-1)}, c - \sum_{\ell=1}^{N^*_n - 1} W_{(\ell, n-1)}\big\} \leq W_{ (N^*_n, n) },
        $
        so we also have that
        $
        W_i \leq W_{ (N^*_n, n) }.
        $
        On the other hand, if $W_{ (N^*_n, n) } < W_{ (N^*_n + 1, n) }$ then the bounds in \eqref{eq:maximands-bounds}
        imply that
        $
        h(\cW_i)  < W_{ (N^*_n+1, n) },
        $
        so we obtain from $W_i \leq h(\cW_i)$ that
        $
        W_i \leq  W_{ (N^*_n, n) }.
        $

    \item[\textsc{Case 2:} $\tau^i_{n-1} = N^*_n$.]
        The left inequality of \eqref{eq:order-statistic-skip-one} with $\ell = N^*_n$ tells us that
        we have two sub-cases to consider here:
        \begin{enumerate*}[label=(\roman*)]%Inline enumerate
            \item when $W_{(N^*_n,n)}$ is equal to $W_{ (N^*_n,n-1) }$, and
            \item when $W_{(N^*_n,n)}$ is strictly smaller than $W_{ (N^*_n,n-1) }$.
        \end{enumerate*}
        In the first sub-case, if $\tau^i_{n-1} = N^*_n$ and $W_{(N^*_n,n)} = W_{ (N^*_n,n-1) }$,
        then the first $N^*_n$ order statistics of $\cW$ and of $\cW_i$ agree and
        $c - \sum_{\ell=1}^{N^*_n} W_{ (\ell, n-1) } = c - \sum_{\ell=1}^{N^*_n} W_{ (\ell, n) } < W_{ (N^*_n + 1, n) }$.
        Thus, if $W_{ (N^*_n, n) } = W_{ (N^*_n + 1, n) }$ then
        $h(\cW_i) = \max\{ W_{ (N^*_n, n) }, c - \sum_{\ell=1}^{N^*_n} W_{ (\ell, n) } \} =  W_{ (N^*_n, n) }$,
        and we are done.
        Otherwise, if $W_{ (N^*_n, n) } < W_{ (N^*_n + 1, n) }$ then
        $h(\cW_i) < W_{ (N^*_n+1, n) }$
        so that $W_i \leq h(\cW_i) < W_{ (N^*_n+1, n) }$ implies that $W_i \leq W_{ (N^*_n, n) }$.
        In the second sub-case, if $\tau^i_{n-1} = N^*_n$ and $W_{(N^*_n,n)} < W_{ (N^*_n,n-1) }$
        then we have that $W_i = W_{ (N^*_n, n) }$, and the result follows.

\end{description}

\emph{Only If.} We now suppose that $W_i \leq W_{ (N^*_n, n) }$, and we show that
$W_i \leq h(\cW_i) \equiv \max\big\{ W_{(\tau^i_{n-1}, n-1)}, c - \sum_{\ell=1}^{\tau^i_{n-1}} W_{(\ell, n-1)}\big\}$
by proving that $W_{ (N^*_n, n) } \leq  h(\cW_i)$.
Just as before, we consider separately the two possible realizations of $\tau^i_{n-1}$.

\begin{description}
    \item[\textsc{Case 1:} $\tau^i_{n-1} = N^*_n - 1$.]
        We have two sub-cases to consider here.
        First, if $W_{ (N^*_n, n) } \leq  W_{ (N^*_n-1, n-1) }$ then the lower bound $W_{ (N^*_n, n) } \leq  h(\cW_i)$ is trivial.
        Second, if $ W_{ (N^*_n-1, n-1) } < W_{ (N^*_n, n) }$ we show that the right maximand is bounded below by $W_{ (N^*_n, n) }$.
        In this instance, the first $N^*_n-1$ order statistic of $\cW$ and $\cW_i$ agree so
        the definition of $N^*_n$ gives us that
        $W_{ (N^*_n,n) } \leq c - \sum_{\ell=1}^{N^*_n-1} W_{(\ell, n)} = c - \sum_{\ell=1}^{N^*_n-1} W_{(\ell, n-1)}$,
        and we are done.

    \item[\textsc{Case 2:} $\tau^i_{n-1} = N^*_n$.]
        If $\tau^i_{n-1} = N^*_n$ the left inequality of \eqref{eq:order-statistic-skip-one} tells us that
        $W_{ (N^*_n, n) } \leq W_{ (N^*_n, n-1) }$, so the lower bound $W_{ (N^*_n, n) } \leq  h(\cW_i)$ immediately follows.
        \Halmos
\end{description}
\endproof\vspace{0.75cm}

The representation \eqref{eq:Nstar-threhsold-representation} for $R^*_n(c,p,r)$
provides us with an easy way for proving that $\E{R^*_n(c,p,r)} \leq npr \F{ \epsilon_{n p}(c) }$.
We just need to note that the expected total reward collected by the offline-sort algorithm
is bounded above by the solution of some appropriate optimization problem.
Our argument does not require independence of item weights.
The threshold equivalence of Lemma \ref{lm:threshold-policy-equivalence-offline-sort}
holds on every sample path, and the relaxation that follows only uses properties of the
weight distribution $F$ and of the arrival probability $p$ \citep[see also][]{Steele:MATHAPPLWARSAW2016}.

\begin{proposition}[Prophet upper bound]\label{prop:prophet-upper-bound}
Consider a knapsack problem with capacity $0 \leq c < \infty$ and
with items that arrive over $1 \leq n <\infty$ periods
according to a Bernoulli process with arrival probability $p\in(0,1]$.
If the items have rewards equal to $r$ and weights with continuous distribution $F$,
then for $\epsilon_{n p}(c) = \sup\left\{ \epsilon \in [0, \infty) :~ \int_{0}^{\epsilon} w \, dF(w) \leq \frac{c}{n p} \right\}$ we have that
\begin{equation}\label{eq:prophet-upper-bound}
\E{ R^*_n (c,p,r) } \leq npr \F{ \epsilon_{n p}(c) }.
\end{equation}
\end{proposition}

\proof{Proof.}
To prove inequality \eqref{eq:prophet-upper-bound},
we begin with two easy cases.
If $c = 0$ then $R^*_n(0,p,r) = 0$, and the bound \eqref{eq:prophet-upper-bound} is trivial.
Similarly, if $\mu = \E{W_1 \mid B_1=1} = \int_{0}^{\infty} w \, dF(w)$ and $n p\mu < c < \infty$
then the definition of the function $\epsilon_{n p}(c)$ tells us that
$\epsilon_{n p}(c) = + \infty$ so $\F{ \epsilon_{n p}(c) }=1$ and the bound \eqref{eq:prophet-upper-bound}
is again trivial because $R^*_n(c,p,r) \leq \sum_{i=1}^n r \1{W_i < \infty}$
%\red{$\E{N^*_n} \leq \E{\sum_{i=1}^n \1{W_i < \infty} } = n p$}
for all $c \in [0, \infty)$,
and this last right-hand side has expected value equal to $npr$.

Next, we consider the case in which $0< c \leq n p\mu$.
If $\cW_i \equiv \{W_1, \ldots, W_{i-1}, W_{i+1}, \ldots, W_n\}$ and $\cG_i = \sigma\{\cW_i\}$ is the $\sigma$-field
generated by the sample $\cW_i$, then we obtain from Lemma \ref{lm:threshold-policy-equivalence-offline-sort}
and from the definition \eqref{eq:offline-sort}
that for each $i\in [n]$ there is a $\cG_i$-measurable threshold $h(\cW_i)$
such that one has the representation as well as the capacity constraint
$$
 R^*_n(c,p,r) = \sum_{i=1}^{n} r\1{ W_i \leq h(\cW_i) }
 \quad \quad \text{and} \quad \quad
 \sum_{i=1}^{n} W_i \1{ W_i \leq h(\cW_i) } \leq c.
$$
In turn, we can obtain an upper bound for $\E{R^*_n(c,p,r)}$ by maximizing
the sum $\sum_{i=1}^n \E{ r\1{ W_i \leq h_{i} } }$ over all thresholds $(h_1, h_2, \ldots, h_n)$
that satisfy an analogous capacity constraint and that have the same measurability property.
Formally, we have the inequality
\begin{eqnarray}\label{eq:optimization-problem}
   \E{ R^*_n(c,p,r) } \leq  & \displaystyle\max_{ (h_1, \ldots, h_n) }   & \sum_{i=1}^n \E{ r\1{ W_i \leq h_{i} } } \\
                           	&  \text{ s.t. }                             & \sum_{i=1}^{n} W_i \1{ W_i \leq h_{i} } \leq c  \quad\text{almost surely}\notag\\
                            &                                            & h_{i} \in \cG_{i} \quad\text{for all } i \in [n]. \notag
\end{eqnarray}
Since $\epsilon_{n p}(c) > 0$ and because the capacity
constraint holds almost surely (and thus also in expectation),
we have the further upper bound
\begin{eqnarray}\label{eq:optimization-problem2}
\E{ R^*_n (c,p,r) }
\leq  &  \displaystyle\max_{ (h_1,\ldots,h_n) }  & \sum_{i=1}^n \E{ r\1{ W_i \leq h_{i} } \{ 1 - \epsilon^{-1}_{n p}(c)W_i\} } + cr \epsilon^{-1}_{n p}(c) \\
      &  \text{ s.t. }              & \sum_{i=1}^{n} W_i \1{ W_i \leq h_{i} } \leq c  \quad\text{almost surely}\notag\\
      &                             & h_{i} \in \cG_{i} \quad\text{for all } i \in [n]. \notag
\end{eqnarray}
Because $h_{i}$ is $\cG_{i}$-measurable, an application of the tower property gives us that
$$
\E{ \E{ r\1{ W_i \leq h_{i} } \{ 1 - \epsilon^{-1}_{n p}(c)W_i\} \mid \cG_i } }
    = pr\E{ \int_0^{ h_{i}} \{ 1 - \epsilon^{-1}_{n p}(c) w \} \, dF(w) },
$$
so, after we drop the two constraints in \eqref{eq:optimization-problem2} we obtain that
\begin{equation}\label{eq:optimization-problem3}
\E{ R^*_n (c,p,r) }
\leq   \mathfrak{p}^* = \max_{ (h_1,\ldots,h_n)  } \sum_{i=1}^n  pr \E{ \int_0^{ h_{i}} \{ 1 - \epsilon^{-1}_{n p}(c) w \} \, dF(w) } + cr \epsilon^{-1}_{n p}(c).
\end{equation}
The maximization problem on the right hand side is separable,
and the quantity $\E{ \int_0^{ h_{i}} \{ 1 - \epsilon^{-1}_{n p}(c) w \} \, dF(w) }$
is maximized by setting $h_{i} = \epsilon_{n p}(c)$ almost surely and for all $i \in [n]$.
Thus, it follows that
\begin{align*}
\mathfrak{p}^*  & = \sum_{i=1}^n   \max_{ h_{i} } pr \E{ \int_0^{ h_{i}} \{ 1 - \epsilon^{-1}_{n p}(c) w \} \, dF(w) } + cr \epsilon^{-1}_{n p}(c)\\
                & = npr \bigg\{ \F{ \epsilon_{n p}(c) } - \epsilon^{-1}_{n p}(c) \bigg[ \int_0^{\epsilon_{n p}(c)} w \, dF(w) - \frac{c}{n p} \bigg] \bigg\}.
\end{align*}
The integral representation \eqref{eq:epsilon-k-integral}
then tells us that the second summand is equal to zero,
so after we recall \eqref{eq:optimization-problem3} we obtain that
$$
\E{ R^*_n (c,p,r) }
\leq \mathfrak{p}^* = npr \F{ \epsilon_{n p}(c) } \quad \quad \text{for all } 0 < c  \leq n p \mu,
$$
completing the proof of \eqref{eq:prophet-upper-bound}.
\halmos\endproof%\vspace{0.75cm}

\section[The reoptimized policy and its value function]
{The reoptimized policy $\boldsymbol{\widehat{\pi}}$ and its value function}\label{se:heuristic-policy}

In the course of proving Proposition \ref{prop:prophet-upper-bound}, we observed that
if $\cG_i = \sigma\{W_1, \ldots, W_{i-1}, W_{i+1}, \ldots, W_k\}$ is the $\sigma$-field
generated by the sample $\{W_1, \ldots, W_{i-1}, W_{i+1}, \ldots, W_k\}$,
then the expected value of the offline solution $R^*_k(x,p,r)$ satisfies the upper bound
\begin{eqnarray*}%\label{eq:optimization-problem}
   \E{ R^*_k(x,p,r) } \leq  & \displaystyle\max_{ (h_1, \ldots, h_k) }   & \sum_{i=1}^k \E{ r\1{ W_i \leq h_{i} } } \\
                           	&  \text{ s.t. }                             & \sum_{i=1}^{k} W_i \1{ W_i \leq h_{i} } \leq x  \quad\text{almost surely}\notag\\
                            &                                            & h_{i} \in \cG_{i} \quad\text{for all } i \in [k]. \notag
\end{eqnarray*}
We also noticed that the optimization problem on the right-hand side can be relaxed by first adding to its objective the
quantity $\epsilon^{-1}_{k p}(x) r \left\{ x - \E{\sum_{i=1}^{k} W_i \1{ W_i \leq h_{i} }} \right\} \geq 0$,
and then by dropping the two constraints.
This then gives us the further upper bound
\begin{equation}\label{eq:optimization-problem3-repeated}
\E{ R^*_k (x,p,r) }
\leq   \max_{ (h_1,\ldots,h_k)  } \sum_{i=1}^k  p r \E{ \int_0^{ h_{i}} \{ 1 - \epsilon^{-1}_{k p}(c) w \} \, dF(w) } + x r \epsilon^{-1}_{k p}(x),
\end{equation}
which is maximized by setting $h_i = \epsilon_{k p}(x)$ for all $i \in [k]$.
We can now use this reoptimized solution for all $x \in [0,\infty)$ and all $1 \leq k  < \infty $
to construct the online feasible threshold policy $\widehat \pi \in \Pi(n,c,p)$.
Specifically, since $\epsilon_{k p}(x)$ may exceed $x$, we set for $p \in (0,1]$
\begin{equation}\label{eq:heuristic-threshold}
\widehat h_{k} (x) = \min\{ x, \epsilon_{k p}(x)\},
\end{equation}
and we define the reoptimized policy $\widehat \pi$ through the threshold $\{ \widehat h_n, \widehat h_{n-1}, \ldots, \widehat h_1\}$.
Thus, if the remaining capacity is $x$ when item $i$ is first presented, then
item $i$ is selected if and only if its weight $W_i \leq \widehat h_{n-i+1}(x)$.

In turn, the threshold functions $\{ \widehat h_k: 1 \leq k <\infty\}$
induce a sequence of value functions $\{ \widehat v_k: [0, \infty) \rightarrow \R_+: 0 \leq k < \infty\}$
such that $\widehat v_{k}(x)$ represents the expected reward to-go of
the reoptimized policy when there are $k$ remaining periods and the
current level of remaining knapsack capacity is $x$.
If $\widehat v_0(x) =0$ for all $x \in [0,\infty)$,
then the value $\widehat v_k(x)$ is given by the recursion
\begin{eqnarray}\label{eq:v-hat}
  \widehat{v}_k(x)
    &=& p \left( 1- \F{ \widehat{h}_{k}(x)} \right) \widehat{v}_{k-1}(x)
    + p \int_0^{\widehat{h}_{k}(x)} \{ r + \widehat{v}_{k-1}(x-w) \} \, d\F{w} + (1-p) \widehat v_{k-1}(x) 	\nonumber\\
	&=& \left( 1- p\F{ \widehat{h}_{k}(x)} \right) \widehat{v}_{k-1}(x) + p \int_0^{\widehat{h}_{k}(x)} \{ r + \widehat{v}_{k-1}(x-w) \} \, d\F{w}.
\end{eqnarray}
By setting the number of remaining periods to $n$ and the knapsack capacity to $c$,
we find that
$$
\widehat{v}_n(c) = \E{ R^{\widehat \pi}_n (c, p, r) }.
$$
To verify the validity of the recursion \eqref{eq:v-hat},
we condition on what happens in the $k$th-to-last period.
With probability $1-p$ the arriving item has arbitrarily large weight
(equivalently, no item arrives),
the number of the remaining periods decreases to $k-1$
and the level of remaining capacity, $x$, stays the same.
This then yields the term $(1-p)\widehat v_{k-1}(x)$ in the first line of \eqref{eq:v-hat}.
On the other hand, with probability $p$ the arriving item has weight distribution $F$,
and we can further condition on its realization, $w$.
If $w > \widehat{h}_{k}(x)$ then the item is rejected,
the level of remaining capacity does not change, and the number of remaining periods decreases by one.
That is, if the item is rejected, the expected reward to-go is given by $\widehat v_{k-1}(x)$ and,
since rejections happen with probability
$p(1-\F{\widehat h_{k}(x)})$,
we recover the first summand on the top line of \eqref{eq:v-hat}.
On the other hand, if $w \leq \widehat{h}_{k}(x)$ the $k$th-to-last item is included in the knapsack.
Such a decision produces an immediate reward of $r$, and it depletes $w$ units of capacity.
The new remaining capacity then becomes $x-w$, and
the number of remaining periods decreases to $k-1$.
The decision maker's payoff for including this item is then given by $r +  \widehat v_{k-1}(x-w)$
and, by integrating this payoff
against the measure $p\,dF(w)$ for $w \in [0, \widehat{h}_{k}(x)]$,
we find the second summand on the first line of the recursion \eqref{eq:v-hat}.

The reoptimized heuristic $\widehat \pi$ then takes the solution of the offline relaxation \eqref{eq:optimization-problem3-repeated}
and turns it into an online algorithm through the threshold $\widehat h_k$ given in \eqref{eq:heuristic-threshold}.
This direct link provides us with enough tractability to be able to quantify the difference in expected performance
between the reoptimized heuristic and the offline solution and---as a result---to prove the logarithmic regret bound.
Instead, the optimal dynamic programming policy cannot be expressed explicitly
and it lacks of the regularity needed to make any meaningful analytical progress.
However, we note here that both the reoptimized heuristic and the optimal dynamic programming policy
can be computed numerically in polynomial time, and we refer the reader to Section \ref{se:numerical}
for more details on our numerical work.

\section{On the typical class}\label{se:typical-class}

The weight distribution $F$ plays a crucial role in the study of the performance of optimal
and near-optimal item selections for the dynamic and stochastic knapsack problem with equal rewards.
Because the weights are not equal, the remaining capacity process
exhibits substantial randomness, and this may lead to unexpected behavior.
As such, regularity conditions on  the weight distribution $F$ are commonplace
in the related literature.
For instance, \citet{CofFlaWeb:AAP1987} only consider distributions $F$ such that $\F{w} \sim A w^\alpha$ as $w \rightarrow 0$
for some $A, \alpha > 0$, while
\citet{BruRob:AAP1991} expand this class to include all of the weight distributions $F$ such
that $\limsup_{w \rightarrow 0^+} \F{\lambda w}/\F{w} < 1$.
Furthermore, \citet[][Section 5]{PapRajKle:MS1996} show that one must require concavity of $F$
to obtain structural properties such as monotonicity of the optimal threshold functions
and concavity of the optimal value functions.

Here, we consider distributions that belong to the typical class characterized in Definition \ref{def:typical-class}.
As we mentioned earlier, this class is broad enough to include most well-known
non-negative continuous distributions.
Such breadth comes from the fact that Conditions
\eqref{eq:cdf-condition} and \eqref{eq:monotonicity-wf(w)} in Definition \ref{def:typical-class}
must hold only on $(0,\bar w)$
for some $\bar w>0$, and that one has the flexibility of choosing different parameter $\bar w$
for different distribution $F$.
For instance, the uniform distribution $f(w) = \1{w \in (0,1) }$
and the exponential distribution $f(w)=\alpha e^{-\alpha w} \1{w > 0}$ are both typical,
but they require different choices of $\bar w$.
For the uniform distribution, Conditions \eqref{eq:cdf-condition} and \eqref{eq:monotonicity-wf(w)}
hold on all of its support and one can choose $\bar w=1$,
while for the exponential distribution, Condition \eqref{eq:monotonicity-wf(w)}
holds only on $(0,3/\alpha)$ and one can set $\bar w=3/\alpha$.
Similarly, one can check that the truncated normal distribution on $(0,b)$ with density
$f(w) = A \exp\{ -  (w - \upsilon)^2 / (2\varsigma^2) \} \1{w \in (0,b)}$ for $\upsilon \in \R$, $\varsigma > 0$,
and $A$ being the appropriate normalizing constant,
is typical with $\bar w = \min\{ \tfrac{1}{2} ( \upsilon  + \sqrt{\upsilon^2 + 12 \varsigma^2}), b \}$.
The truncated logistic distribution on $(0,b)$ and the logit-normal distribution
are additional examples of typical distributions,
though the respective $\bar w$'s have to do with the smallest positive root of related
transcendental equations.
The families of distributions listed below also belong to the typical class.
\begin{enumerate}
	
	\item   \emph{Power distributions.}
	Distributions such that $\F{w} = A w^{\alpha}$ for some $A,\alpha>0$ on $(0, \bar w)$
	are typical. Condition \eqref{eq:cdf-condition} is immediately verified.
	The function $w^3f(w) = A \alpha w^{\alpha + 2}$ is increasing because $A,\alpha > 0$,
	so \eqref{eq:monotonicity-wf(w)} holds as well.
	
	\item   \emph{Convex distributions.}
	Distributions $F$ that are convex in a neighborhood of 0 and that have continuous density $f$ are typical.
	Convexity tells us that $\F{\lambda w} \leq \F{w} \lambda$ so \eqref{eq:cdf-condition} follows.
	Furthermore, convexity also gives us that the density $f$ is non-decreasing, so \eqref{eq:monotonicity-wf(w)} is verified.
	
	\item   \emph{Mixtures of typical distributions.}
	The class of typical distributions is closed under mixture.
	If $F$ and $G$ are two typical distributions and $\beta\in[0,1]$ then
	it is easy to see that the mixture distribution $\beta F + (1-\beta)G$ is also typical.

\end{enumerate}

It is important to note, however, that one can construct examples of distributions
that do not belong to the typical class.
For instance, the distribution  $\F{ w } =  \frac{ \log \bar w  }{ \log w }$ for $\bar w < 1$ and $w \in (0,\bar w)$
is an example that satisfies Condition \eqref{eq:monotonicity-wf(w)}
but violates Condition \eqref{eq:cdf-condition}.
For a fixed $0 < \lambda < 1$ , one can easily check that
$$
\limsup_{w \rightarrow 0^+} \frac{\F{ \lambda w }}{ \F{ w } } = \limsup_{w \rightarrow 0^+} \frac{\log w}{\log \lambda + \log w} = 1,
$$
so Condition \eqref{eq:cdf-condition} fails to hold.
On the other hand, the function $w^3 f(w) = - \frac{w^2 \log \bar w }{(\log w)^2}$ is increasing on $(0, \bar w)$
and Condition \eqref{eq:monotonicity-wf(w)} is satisfied.

The distribution $\F{ w } = A \int_0^w \{ \sin \left( 1/u \right) \}^2 \, du$ for $w \in (0, \bar w)$ and $A = ( \int_0^{\bar w} \{ \sin \left( 1/u \right) \}^2 \, du)^{-1} > 0$ is an example that satisfies Condition \eqref{eq:cdf-condition} while violating Condition \eqref{eq:monotonicity-wf(w)}.
In fact, one has that the limit
$$
\limsup_{w \rightarrow 0^+} \frac{\F{ \lambda w }}{ \F{ w }} = \lambda < 1,
$$
but the function $w^3 f(w) = A w^3 \{\sin (1/w)\}^2$ oscillates infinitely many times in a (positive) neighborhood of zero,
so the monotonicity \eqref{eq:monotonicity-wf(w)} fails to hold.

We conclude this section by observing that Condition \eqref{eq:cdf-condition}
regarding the behavior of $F$ at zero
is equivalent to the condition required by \citet{BruRob:AAP1991},
and by proving that we can equivalently state it as a property of the ratio $wF(w)/\int_0^w u \ dF(u)$.
This equivalent property will be important to our analysis.

\begin{lemma}[Equivalence of CDF Conditions]\label{lm:cdf-conditions-equivalence}
  There are constants $0 < \lambda < 1$ and $0 < \gamma < 1$
  and a value $\bar{w} > 0$ such that
  \begin{equation}\label{eq:cdf-condition-in-lemma}
    \frac{\F{ \lambda w }}{ \F{  w  }} \leq \gamma < 1 \quad \quad \text{for all } w \in (0, \bar{w})
  \end{equation}
  if and only if there is a constant $1 < M < \infty$ such that
  \begin{equation}\label{eq:cdf-condition-integral-representation}
    \frac{w \F{ w }}{\int_{0}^{w} u  \, d\F{u} } \leq M < \infty \quad \quad \text{for all } w \in (0, \bar{w}).
  \end{equation}
\end{lemma}

\proof{Proof of Lemma \ref{lm:cdf-conditions-equivalence}.}
\emph{If.} Suppose there is a constant $1 < M < \infty$
such that condition \eqref{eq:cdf-condition-integral-representation} holds.
Next, note that for any $\lambda \in (0,1)$ and any $w \in (0, \bar{w})$
one has the bounds
$$
0 \leq \int_{0}^{w} u \, d\F{u}
    \leq \lambda w \int_{0}^{\lambda w}  \, d\F{u}  + w \int_{\lambda w}^{w}  \, d\F{u}
    = w \F{ w } - w \F{ \lambda w }(1 - \lambda),
$$
so it follows that
$$
\frac{\F{ w }}{\F{ w } - \F{ \lambda w }(1 - \lambda)}
\leq \frac{w \F{ w }}{\int_{0}^{w} u  \, d\F{u} }.
$$
In turn, condition \eqref{eq:cdf-condition-integral-representation} tells us that there is $1 < M < \infty$
such that the right-hand side above is bounded by $M$ so, after rearranging, we obtain that
$$
\frac{\F{ \lambda w }}{\F{w}}  \leq \frac{M-1}{M(1-\lambda)}  \quad \quad \text{for all } w \in (0, \bar{w}).
$$
Condition \eqref{eq:cdf-condition-in-lemma} then follows after one chooses any $\lambda < M^{-1}$
and sets $\gamma = (M-1) / [M(1-\lambda)] < 1$.

\medskip
\emph{Only if.} Suppose that there are constants $0< \lambda < 1$
and $0 < \gamma < 1$ such that condition \eqref{eq:cdf-condition-in-lemma} holds for some $\bar{w}>0$.
Then we have that
$$
0 < 1 - \gamma \leq  1 - \frac{ \F{ \lambda w }}{ \F{ w } } = \int_{\lambda w}^{w} \frac{d \F{u} }{ \F{ w }}
\quad \quad \text{for all } w \in (0, \bar{w}).
$$
Moreover, if we multiply both sides by $\lambda w$ and use the fact that $\lambda w \leq u$
for all $u \in (\lambda w, w)$ we also have that
\begin{equation*}%\label{eq:lower-bound-lemma-intermediate}
\lambda w ( 1- \gamma) \leq \int_{\lambda w}^{w} \frac{ \lambda w }{\F{ w }} \, d \F{u}
 \leq \int_{0}^{w} \frac{u }{\F{ w }} \, d \F{u}.
\end{equation*}
Next, we divide both sides by $w$ and rearrange to obtain that
$$
\frac{w \F{ w }}{\int_{0}^{w} u \, d \F{u}} \leq \frac{1}{\lambda (1-\gamma)} \quad \quad \text{for all } w \in (0, \bar{w}),
$$
so condition \eqref{eq:cdf-condition-integral-representation}
follows by setting $M = [\lambda ( 1- \gamma)]^{-1}$,
and the proof is now complete.
\Halmos\endproof%\vspace{0.75cm}

\section{A logarithmic regret bound}\label{se:main-proof}

To prove that the regret grows at most logarithmically, we let
\begin{equation}\label{eq:K}
K = \ceil{ \frac{c}{p \int_{0}^{\bar w} w f(w) \, dw}}
\end{equation}
and focus on dynamic and stochastic knapsack problems with more than $K$ periods.
Of course, this is without loss of generality because the  quantity $K$ defined in \eqref{eq:K}
is a constant that does not depend on the number of periods $n$, so we can ignore the
last $K$ decisions without affecting our regret bound.
When $k \geq K$ we have
\begin{enumerate*}[label=(\roman*)]
    \item that $\epsilon_{k p}(x) \leq \bar w$ for all $x \in [0,c]$, and
    \item that the integral representation \eqref{eq:epsilon-k-integral} always holds.
\end{enumerate*}
Thus, we are focusing on problem instances in which
we can use the properties of the typical class in full.

In our proof, we will repeatedly use the following two properties of the consumption function $\epsilon_{k p}(x)$.
First, we obtain from definition \eqref{eq:epsilon-k-definition} that
the consumption functions are non-increasing in $k$. That
is, for $p \in(0,1]$ one has the monotonicity
\begin{equation}\label{eq:epsilon-k-monotone}
\epsilon_{(k+1)p}(x) \leq \epsilon_{k p}(x)
\quad \quad \text{for all } x \in [0,\infty) \text{ and all } k \geq 1.
\end{equation}
Second, provided that the weight distribution $F$ has continuous density $f$,
an application of the implicit function theorem gives us that the function $\epsilon_{k p}(x)$ is differentiable on $(0, k p\mu)$,
and that its first derivative $\epsilon'_{k p}(x)$ is given by
\begin{equation}\label{eq:epsilon-k-first-derivative}
  \epsilon'_{k p}(x) = \frac{1}{k p \epsilon_{k p}(x) f(\epsilon_{k p}(x))}
  \quad \quad \text{if } 0 < x < k p\mu.
\end{equation}
The proof of the regret bound then comes in two parts.
In the next section we derive several estimates that have to do
with the weight distribution belonging to the typical class and with $k\geq K$,
while in Section \ref{se:analysis-of-residuals} we
estimate the gap $kpr \F{ \epsilon_{k p}(x) } - \widehat{v}_k(x)$.

\subsection{Preliminary observations}

When $k \geq K$ the properties that characterize typical weight distributions
can be used to obtain general estimates that are crucial to our analysis.
As a warm-up we obtain the following estimate
on the mismatch between the probability of an item weight being smaller than the
feasible threshold $\widehat h_{k}$ and the probability of the same weight
being smaller than the consumption function $\epsilon_{k p}$.

\begin{lemma}\label{lm:mismatch-epsilon-hhat}
If the weight distribution $F$ belongs to the typical class then
there is $1 < M < \infty$ depending only on $F$ such that
\begin{equation}\label{eq:key-upper-bound}
\frac{k p \epsilon_{k p}(x) \F{ \epsilon_{k p}(x)} }{x}
 \leq M
 \quad \text{ for all $x \in (0,c]$, $p \in (0,1]$, and all } k \geq K \equiv \ceil{ \frac{c} { p\int_{0}^{\bar w} w f(w) \, dw} }.
\end{equation}
In turn, we also have that
\begin{equation}\label{eq:difference-Fepsiloon-Fh}
  \F{ \epsilon_{k p}(x) } - \F{ \widehat{h}_{k}(x) } \leq \frac{M}{k p}
  \quad \text{ for all $x \in [0,c]$, $p \in (0,1]$, and all } k \geq K.
\end{equation}
\end{lemma}

\proof{Proof.}
The uniform bound \eqref{eq:key-upper-bound} is essentially a restatement of
inequality \eqref{eq:cdf-condition-integral-representation} in Lemma \ref{lm:cdf-conditions-equivalence}.
If $x \in (0, c]$ and $k \geq K$, then we have that
$$
\frac{x}{k p} \leq \frac{x}{Kp} \leq \frac{c}{Kp} \leq \int_{0}^{\bar{w}} w f(w) \, dw \leq \mu,
$$
so the definition \eqref{eq:epsilon-k-definition}
of the consumption function $\epsilon_{k p}(x)$ and the equality \eqref{eq:epsilon-k-integral} give us that
\begin{equation}\label{eq:observations-post-K}
\epsilon_{k p}(x) \leq \bar{w}
\quad \text{and} \quad
\int_{0}^{\epsilon_{k p}(x)} w f(w) \, dw = \frac{x}{k p}
\quad \quad \quad \text{ for all $k \geq K$  and all $x \in (0,c]$.}
\end{equation}
The two observations in \eqref{eq:observations-post-K} together with the bound \eqref{eq:cdf-condition-integral-representation}
in which we replace $w$ with $\epsilon_{k p}(x)$ then imply that
$$
\frac{k p \epsilon_{k p}(x) \F{ \epsilon_{k p}(x)} }{x}
= \frac{\epsilon_{k p}(x) \F{ \epsilon_{k p}(x)} }{\int_{0}^{\epsilon_{k p}(x)} u f(u) \, du}
\leq M
 \quad \text{ for all $k \geq K$  and $x \in (0,c]$,}
$$
concluding the proof of the uniform bound \eqref{eq:key-upper-bound}.

We now turn to inequality \eqref{eq:difference-Fepsiloon-Fh}.
If $x=0$ then inequality \eqref{eq:difference-Fepsiloon-Fh} is obvious.
Otherwise, if $x>0$ we recall from \eqref{eq:heuristic-threshold} that $\widehat{h}_{k}(x) = \min\{ x, \epsilon_{k p}(x)\}$,
so  the left-hand side of \eqref{eq:difference-Fepsiloon-Fh} is equal to 0
when $\epsilon_{k p}(x) \leq x < \infty$, and inequality \eqref{eq:difference-Fepsiloon-Fh} is again trivial.
Instead, if $0< x < \epsilon_{k p}(x)$, we obtain from \eqref{eq:key-upper-bound}
that
$$
\F{ \epsilon_{k p}(x) } \leq \frac{M x}{k p \epsilon_{k p}(x)} \leq \frac{M}{k p}
\quad \text{ for all $k \geq K$  and $0 < x < \epsilon_{k p}(x)$,}
$$
concluding the proof of the lemma.
\Halmos\endproof\vspace{0.75cm}

In the same spirit of Lemma \ref{lm:mismatch-epsilon-hhat}, we can
also estimate the difference in the probability of selecting an upcoming
item as a function of the number of remaining periods.

\begin{lemma}\label{lm:F-epsilon-difference-bound}
  For $p\in(0,1]$, all $x \in [0,c]$, and all $k \geq K$ we have that
  $$
  \F{ \epsilon_{(k+1)p}(x) } - \F{ \epsilon_{k p}(x) } \leq - \frac{x}{k(k+1)p \epsilon_{k p}(x)}.
  $$
\end{lemma}

\proof{Proof.}
  For $k \geq K$ the equality  \eqref{eq:epsilon-k-integral} and
  the monotonicity \eqref{eq:epsilon-k-monotone}
  give us the representation
  $$
  \frac{x}{k p} - \frac{x}{(k+1)p} = \int_{\epsilon_{(k+1)p}(x)}^{\epsilon_{k p}(x)} wf(w)  \, dw,
  \quad \quad \text{for all } x \in [0,c].
  $$
  If we now replace the integrand $w f(w)$ with the upper bound $\epsilon_{k p}(x)f(w)$
  and rearrange, we obtain
  $$
  \frac{x}{k(k+1)p} \leq \epsilon_{k p}(x) \left[ \F{ \epsilon_{k p}(x) } - \F{ \epsilon_{(k+1)p}(x) } \right],
  $$
  completing the proof of the lemma.
\Halmos\endproof\vspace{0.75cm}

Typical weight distributions are also nice because
one can tightly approximate the difference $\F{ \epsilon_{k p}(x) } - \F{ \epsilon_{k p}(x-w)}$
that accounts for the sensitivity in the remaining capacity of the probability of selecting
the $k$th-to-last item.
A formal estimate is given in the next proposition, and it constitutes a key step in our argument.

\begin{proposition}\label{pr:upper-bound-ratio-cdfs}
    If $p\in(0,1]$ and if the weight distribution $F$ belongs to the typical class,
    then there is a constant $1 < M < \infty$ depending only on $F$ such that one
    has the inequality
    \begin{equation}\label{eq:upper-bound-ratio-cdfs}
      1 - \frac{\F{ \epsilon_{k p}(x - w) }}{ \F{ \epsilon_{k p}(x)} }
      \leq \frac{w^2}{x^2} (1 - M^{-1}) + \frac{w}{k p \epsilon_{k p}(x) \F{ \epsilon_{k p}(x) } }
    \end{equation}
    for all $w \in [0,x]$, $x \in (0,c]$, and all $k \geq K \equiv \ceil{ \frac{c}{ p\int_{0}^{\bar w} w f(w) \, dw} }$.
\end{proposition}

The proof of Proposition \ref{pr:upper-bound-ratio-cdfs} requires
the following intermediate estimate.

\begin{lemma}[Convexity upper bound]\label{lm:cdf-difference-integral-representation}
If $p\in(0,1]$ and if the weight distribution $F$ has continuous density $f$ then
for all $k \geq K$, $x \in [0,c]$ and $ y \in [0,1]$
we have the integral representation
\begin{equation}\label{eq:integral-representation}
  k p \F{  \epsilon_{k p}(x) } - k p \F{  \epsilon_{k p}(x(1- y)) }
  = \int_{x(1- y)}^{x} \frac{1}{\epsilon_{k p}(u)} \, du.
\end{equation}
Moreover, if the distribution $F$ belongs to the typical class
the map $x \mapsto \epsilon_{k p}(x)^{-1}$ is convex on $(0,c)$,
so we also have the upper bound
\begin{equation}\label{eq:cdf-difference-upper-bound}
  k p \F{ \epsilon_{k p}(x) } - k p \F{  \epsilon_{k p}(x(1- y)) }
  \leq \frac{x  y}{2}\left[\frac{1}{\epsilon_{k p}(x)} + \frac{1}{\epsilon_{k p}(x(1- y))} \right].
\end{equation}
\end{lemma}

\proof{Proof.}
Since the weight distribution $F$ has continuous density and $\frac{c}{p\mu} \leq K \leq k$,
we recall from \eqref{eq:epsilon-k-first-derivative} the first derivative
$$
\epsilon'_{k p}(x) = \frac{1}{k p \epsilon_{k p}(x) f(\epsilon_{k p}(x))} \quad \quad \text{for all $x \in (0,c)$}.
$$
The map $x \mapsto \F{ \epsilon_{k p}(x) }$ is then differentiable on $(0,c)$,
and one has that
$$
  \left( k p \F{ \epsilon_{k p}(x) } \right)' = k p \epsilon'_{k p}(x) f(\epsilon_{k p}(x)) = \frac{1}{\epsilon_{k p}(x)}
  \quad \quad \text{for all $x \in (0,c)$}.
$$
In turn, the fundamental theorem of calculus
tells us that for $ y \in [0,1]$ we have the integral
representation
\begin{equation*}
   k p \F{  \epsilon_{k p}(x)  }  - k p  \F{  \epsilon_{k p}(x(1- y)) }
   = \int_{x(1- y)}^{x} \frac{1}{\epsilon_{k p}(u)} \, du,
\end{equation*}
proving the first assertion of the lemma.

To check the convexity of the map $x \mapsto \epsilon_{k p}(x)^{-1}$ ,
we use the expression of the first derivative \eqref{eq:epsilon-k-first-derivative} one more time
to obtain for $k \geq K$ that
\begin{equation*}
  \left( \frac{1}{\epsilon_{k p}(x)}\right)'
  = - \frac{\epsilon'_{k p}(x)}{\epsilon^2_{k p}(x)}
  = - \frac{1}{k p \epsilon^3_{k p}(x) f(\epsilon_{k p}(x))}.
\end{equation*}
If $F$ belongs to the typical class and $k\geq K$
then the monotonicity condition \eqref{eq:monotonicity-wf(w)}
implies that the first derivative $( 1/ \epsilon_{k p}(x))'$
is non-decreasing on $(0,c)$,
so the map $x \mapsto \epsilon_{k p}(x)^{-1}$ is convex.
This convexity property then provides us with a linear majorant
$$
m_{k p}(u) = \frac{u-x}{ y x} \left( \frac{1}{\epsilon_{k p}(x)} - \frac{1}{\epsilon_{k p}((1- y)x)}\right) +  \frac{1}{\epsilon_{k p}(x)}
$$
such that
$$
\frac{1}{\epsilon_{k p}(u)} \leq m_{k p}(u)
\quad \quad \text{for all } u \in [(1- y)x, x].
$$
The representation \eqref{eq:integral-representation} and the integration of the majorant $m_{k p}(u)$ over $[(1- y)x, x]$ give us the upper bound \eqref{eq:cdf-difference-upper-bound}, and the proof of the lemma follows.
\Halmos\endproof\vspace{0.75cm}

We now have all of the estimates we need to complete the proof of Proposition \ref{pr:upper-bound-ratio-cdfs}.

\proof{Proof of Proposition \ref{pr:upper-bound-ratio-cdfs}.}
If $w = 0$ then inequality \eqref{eq:upper-bound-ratio-cdfs} is trivial.
Otherwise, for $K \leq k < \infty$ we consider the function
$g_k : (0, c] \times (0,1] \rightarrow \R$
given by
$$
g_k(x, y)
= \frac{1}{y^2}
\left\{ 1 - \frac{ \F{ \epsilon_{k p}\left( x ( 1 - y) \right) } }{ \F{  \epsilon_{k p}(x) } } - \frac{x y}{k p \epsilon_{k p}(x) \F{ \epsilon_{k p}(x) }} \right\},
$$
and we note that inequality \eqref{eq:upper-bound-ratio-cdfs}
follows by setting $y = w/x \leq 1$ and rearranging,
provided that one has the uniform bound
\begin{equation}\label{eq:g-k-upper-bound}
  g_k(x,y) \leq 1 - M^{-1}
  \quad \quad \text{for all $x \in (0, c]$, $y \in (0,1]$, and $k \geq K$.}
\end{equation}
The function $g_k(x,y)$ is differentiable with respect to $y$ for any given $x \in (0,c]$,
and the $y$-derivative of $g_k(x, y)$ can be written as
$$
\frac{\partial}{\partial y} g_k(x,y)
= \frac{2}{ y^3  k p \F{ \epsilon_{k p}(x) } }
\left\{ \frac{x y}{2} \left[ \frac{1}{\epsilon_{k p}(x)} + \frac{1}{\epsilon_{k p}(x(1-y))}\right]
- k p \F{ \epsilon_{k p}(x) } + k p \F{ \epsilon_{k p}(x(1-y)) } \right\}.
$$
Since $\frac{2}{ y^3  k p \F{ \epsilon_{k p}(x) } } \geq 0$, inequality \eqref{eq:cdf-difference-upper-bound} of Lemma \ref{lm:cdf-difference-integral-representation}
then tells us that the $y$-derivative of $g_k(x,y)$
is non-negative so that the map $y \mapsto g_k(x, y)$ is non-decreasing in $y$
for any given $x \in (0,c]$.
In turn, we have that
$$
g_k(x, y) \leq g_k(x, 1) = 1 - \frac{x}{k p \epsilon_{k p}(x) \F{ \epsilon_{k p}(x) } },
$$
so inequality \eqref{eq:g-k-upper-bound} follows
from the uniform bound \eqref{eq:key-upper-bound},
and the proof of the proposition is now complete.
\Halmos\endproof%\vspace{0.75cm}

\subsection{Analysis of residuals}\label{se:analysis-of-residuals}

To estimate the gap between the expected total reward
collected by the reoptimized policy $\widehat \pi \in \Pi(n,c,p)$ and the
prophet upper bound $npr \F{ \epsilon_{n p}(c)}$,
we study appropriate residual functions.
Specifically, we let
\begin{equation}\label{eq:residual}
  \rho_k(x) = kpr \F{ \epsilon_{k p}(x) } - \widehat{v}_k(x)
  \quad \quad \text{for } x \in [0,c] \text{ and } 1 \leq k \leq n
\end{equation}
be the \emph{residual function} when there are $k$ remaining periods
and the level of remaining capacity is $x$.
The residual function $\rho_k(x)$ is continuous and defined on a compact interval,
so if we maximize with respect to $x$ we obtain the \emph{maximal residual}
\begin{equation}\label{eq:max-residual}
  \widebar \rho_k =  \max_{0 \leq x \leq c } \rho_k(x)
  \quad \quad \text{for }  k \in [n].
\end{equation}

The second half of Theorem \ref{thm:main} is
just a corollary of the following proposition, which verifies that the
maximal residual $\widebar \rho_n = O(\log n)$
as $n \rightarrow \infty$.

\begin{proposition}\label{pr:max-residual-bound}
  If the weight distribution $F$ belongs to the typical class, then there is a constant $1 < M < \infty$
  depending only on the distribution $F$, the arrival probability $p$, and the reward $r$
  such that the maximal residual
  $$
  \widebar \rho_n = \max_{0 \leq x \leq c }\{ npr \F{ \epsilon_{n p}(x) } - \widehat{v}_n(x) \} \leq M + M \log n
  \quad \quad \text{for all } n \geq 1.
  $$
\end{proposition}

For the proof of this proposition we write the maximal residual $\widebar \rho_n$
as a telescoping sum, and we obtain an appropriate upper bound for each summand.
The upper bound follows from the following lemma.

\begin{lemma}\label{lm:residual-time-difference-bound}
  If the weight distribution $F$ belongs to the typical class,
  then there is a constant $1 < M < \infty$
  that depends only on $F$ and the reward $r$ such that the difference
  $$
  \rho_{k+1} (x) - \widebar \rho_k \leq \frac{M}{k+1}
  \quad \quad \text{ for all } x \in [0,c] \text{ and all } k \geq K.
  $$
\end{lemma}

\proof{Proof.}
  The residual function $\rho_k(x)$ defined in \eqref{eq:residual}
  provides us with an alternative representation for the value function $\widehat{v}_{k+1}(x)$
  which gives us the expected total reward selected by policy $\widehat \pi$
  with $k+1$ periods remaining and current knapsack capacity $x$.
  Specifically, if we substitute $\widehat v_k(x)$ with $kpr \F{ \epsilon_{k p}(x) } - \rho_k(x)$
  in the recursion \eqref{eq:v-hat},
  we then obtain that
  \begin{align*}
  \widehat{v}_{k+1}(x)
    =&  \{ 1- p\F{ \widehat{h}_{k+1}(x)} \} \{kpr \F{ \epsilon_{k p}(x) } - \rho_k(x)\}\\
     &  + p\int_0^{\widehat{h}_{k+1}(x)} \{ r + kpr \F{ \epsilon_{k p}(x-w) }  - \rho_k(x-w) \} f(w) \, dw.
  \end{align*}
  Next, if we replace the residuals $\rho_k(\cdot)$ with their maximal value $\widebar \rho_k$ and rearrange,
  we obtain the lower bound
  \begin{eqnarray}\label{eq:v+r}
    kpr \F{ \epsilon_{k p}(x) } &+& pr\F{ \widehat{h}_{k+1}(x) } 	
    + k p^2r \int_0^{\widehat{h}_{k+1}(x)} \{ \F{ \epsilon_{k p}(x-w) } - \F{ \epsilon_{k p}(x) } \} f(w) \, dw 	\\
    & \leq & \widehat v_{k+1}(x) + \widebar \rho_k.		\nonumber
  \end{eqnarray}
  In turn, the definition \eqref{eq:residual} of the residual function tells us that
  $$
  \rho_{k+1}(x) - \widebar \rho_k = (k+1)pr \F{ \epsilon_{(k+1)p}(x) } - (\widehat v_{k+1}(x) + \widebar \rho_k),
  $$
  so if we replace the sum $\widehat v_{k+1}(x) + \widebar \rho_k$ with its lower bound \eqref{eq:v+r}
  and rearrange, we obtain the upper bound
  \begin{eqnarray}\label{eq:residual-time-difference}
    \rho_{k+1}(x) - \widebar \rho_k
    & \leq & pr\left\{ (k+1) \F{ \epsilon_{(k+1)p}(x) } - k \F{ \epsilon_{k p}(x) } - \F{ \widehat{h}_{k+1}(x) } \right\}\\
         && + k p^2r \F{ \epsilon_{k p}(x) } \int_0^{\widehat{h}_{k+1}(x)} \bigg\{ 1 - \frac{ \F{ \epsilon_{k p}(x-w)}  }{\F{ \epsilon_{k p}(x)} } \bigg\} f(w) \, dw. \notag
  \end{eqnarray}
  Next, we obtain from \eqref{eq:upper-bound-ratio-cdfs} that the integral that appears on the right-hand side of \eqref{eq:residual-time-difference}
  satisfies the upper bound
  $$
  \cI_k(x) \!\equiv\!  \int_0^{\widehat{h}_{k+1}(x)} \!\bigg\{\! 1 - \frac{ \F{ \epsilon_{k p}(x-w) } }{\F{ \epsilon_{k p}(x) } } \!\bigg\}\! f(w) dw
           \!\leq\!    \int_0^{\widehat{h}_{k+1}(x)} \!\bigg\{\! \frac{w^2}{x^2}(1-M^{-1}) + \frac{w}{k p \epsilon_{k p}(x) \F{ \epsilon_{k p}(x) } } \!\bigg\}\! f(w) dw.
  $$
  For $w \in [0, \widehat{h}_{k+1}(x)]$ we have the trivial bound $w^2 \leq w \widehat{h}_{k+1}(x)$
  so if we replace $w^2$ with its upper bound $w \widehat{h}_{k+1}(x)$
  on the right-hand side above and integrate
  we obtain that there is $1 < M < \infty$ depending only on $F$ such that
  $$
  \cI_k(x)
  \leq \left[ (1- M^{-1}) \frac{\widehat{h}_{k+1}(x)}{x^2} + \frac{1}{k p \epsilon_{k p}(x) \F{\epsilon_{k p}(x)}} \right] \int_{0}^{\widehat{h}_{k+1}(x)} w f(w) \, dw.
  $$
  We now multiply both sides by $k p^2r \F{ \epsilon_{k p}(x) }$ and simplify to obtain that
  $$
  k p^2r \F{ \epsilon_{k p}(x) } \cI_k(x) \leq \left[ k p^2r \F{ \epsilon_{k p}(x)} (1- M^{-1}) \frac{\widehat{h}_{k+1}(x)}{x^2} + \frac{pr}{\epsilon_{k p}(x)}\right] \int_{0}^{\widehat{h}_{k+1}(x)} w f(w) \, dw.
  $$
  The definition of $\widehat{h}_{k+1}(x) = \min\{ x, \epsilon_{(k+1)p}(x)\}$
  and the monotonicity \eqref{eq:epsilon-k-monotone}
  tell us that $\widehat h_{k+1}(x) \leq \epsilon_{(k+1)p}(x) \leq \epsilon_{k p}(x)$,
  so we obtain a further upper bound
  if we replace the first $\widehat{h}_{k+1}(x)$ on the last right-hand side with $\epsilon_{k p}(x)$
  and the second one with $\epsilon_{(k+1)p}(x)$.
  When we perform these replacements and recall the equality \eqref{eq:epsilon-k-integral},
  we find that
  $$
  k p^2r \F{ \epsilon_{k p}(x) } \cI_k(x)
  \leq \frac{r(1- M^{-1})}{k+1} \frac{k p \epsilon_{k p}(x) \F{ \epsilon_{k p}(x) } }{x} +  \frac{rx}{(k+1) \epsilon_{k p}(x)}.
  $$
  If we now apply the uniform upper bound \eqref{eq:key-upper-bound}
  to the first summand on the right-hand side,
  and rearrange, we obtain that
  $$
  k p^2r \F{ \epsilon_{k p}(x) }  \cI_k(x) \leq \frac{r(M-1)}{k+1} + \frac{rx}{(k+1) \epsilon_{k p}(x)}.
  $$
  We now replace the last summand of \eqref{eq:residual-time-difference}
  with the upper bound above and rearrange to obtain that
  \begin{align*}%\label{eq:residual-time-difference2}
    \rho_{k+1}(x) - \widebar \rho_k
    \leq&  \frac{r(M-1)}{k+1} + kpr \left\{ \F{ \epsilon_{(k+1)p}(x) } - \F{ \epsilon_{k p}(x) } +\frac{x}{k(k+1)p \epsilon_{k p}(x)}\right\} 	\\
         & + pr\left\{ \F{ \epsilon_{(k+1)p}(x)} - \F{ \widehat{h}_{k+1}(x)} \right\}.
  \end{align*}
  Here, Lemma \ref{lm:F-epsilon-difference-bound} tells us that the second summand on the right-hand side
  is non-positive, and inequality \eqref{eq:difference-Fepsiloon-Fh} tells us
  that there is $1 < M < \infty$ depending only on $F$ such that the
  difference $\F{ \epsilon_{(k+1)p}(x) } - \F{ \widehat{h}_{k+1}(x) } $ is bounded
  above by $M/((k+1)p)$.
  When we assemble these observations, we finally find that
  $$
  \rho_{k+1}(x) - \widebar \rho_k \leq \frac{(2M - 1)r}{k+1}
  \quad \quad \text{ for all } x \in [0,c] \text{ and all } k \geq K,
  $$
  concluding the proof of the lemma.
\Halmos\endproof\vspace{0.75cm}

We now have all of the tools we need to complete the proof of Proposition \ref{pr:max-residual-bound} that follows next.

\proof{Proof of Proposition \ref{pr:max-residual-bound}.}
  We write the maximal residual $\widebar \rho_n$ in \eqref{eq:max-residual}  as a telescoping sum and use the definition \eqref{eq:residual}
  of the residual function to obtain that
  $$
  \widebar \rho_n = \widebar \rho_{K} + \sum_{k = K}^{n-1} \{ \widebar \rho_{k+1} - \widebar \rho_k \}
  \leq K +  \sum_{k = K}^{n-1} \{ \widebar \rho_{k+1} - \widebar \rho_k \}.
  $$
  Lemma \ref{lm:residual-time-difference-bound} then tells us  that
  $$
  \widebar \rho_{k+1} - \widebar \rho_k \leq \frac{M}{k+1}
  \quad \quad \text{for all } K \leq k \leq n,
  $$
  so when we combine the last two observations we obtain
  that there is a constant $1 < M < \infty$ that depends only on $F$, $p$, and $r$ such that
  $$
  \widebar \rho_n \leq M + M \log n,
  $$
  just as needed.
\Halmos\endproof%\vspace{0.75cm}

\section{Numerical experiments}\label{se:numerical}

\begin{figure}[t!]
	
	\caption{\textbf{Gap between the prophet upper bound and offline sort for three weight distributions.}}\label{fig:simulations}
	
	\begin{center}
		\includegraphics[width=0.5\textwidth]{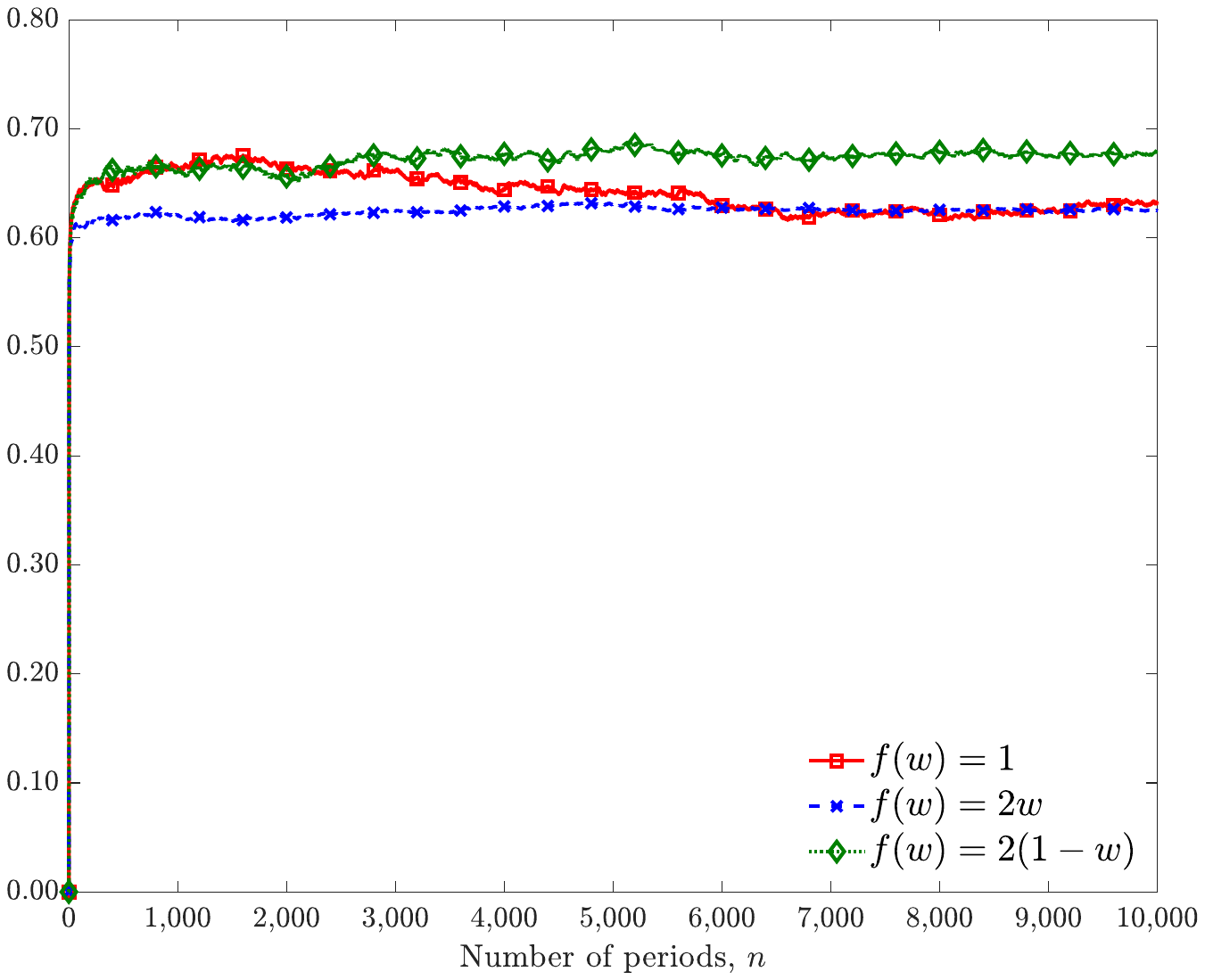}	
        \\[3pt]
	\end{center}	
	
	\scriptsize{\emph{Notes.}
		Difference between the prophet upper bound, $nF(\epsilon_n(1))$,
		and the simulated average (with $100,000$ trials) of the offline solution, $R^*_n(1,1,1)$,
		for three different distributions on the unit interval:
		$f(w) = \1{w \in (0,1) }$, $f(w)=2w\1{w\in (0,1) }$ and $f(w)=2(1-w)\1{w\in (0,1) }$.
		In each case we take the arrival probability $p=1$, the knapsack capacity $c=1$,
		the reward $r=1$, and we vary the number of periods $n \in \{1,2,\ldots,10000\}$.
		The chart suggests that the gap between the prophet upper bound and the simulated average
		of the offline solution does not grow with $n$.
	}	
\end{figure}

\begin{figure}[th!]
	
	\caption{\textbf{Value functions and scaled regret bounds for three weight distributions}}\label{fig:value-function-and-scaled-regret}
	
	\begin{center}
		\includegraphics[width=0.485\textwidth]{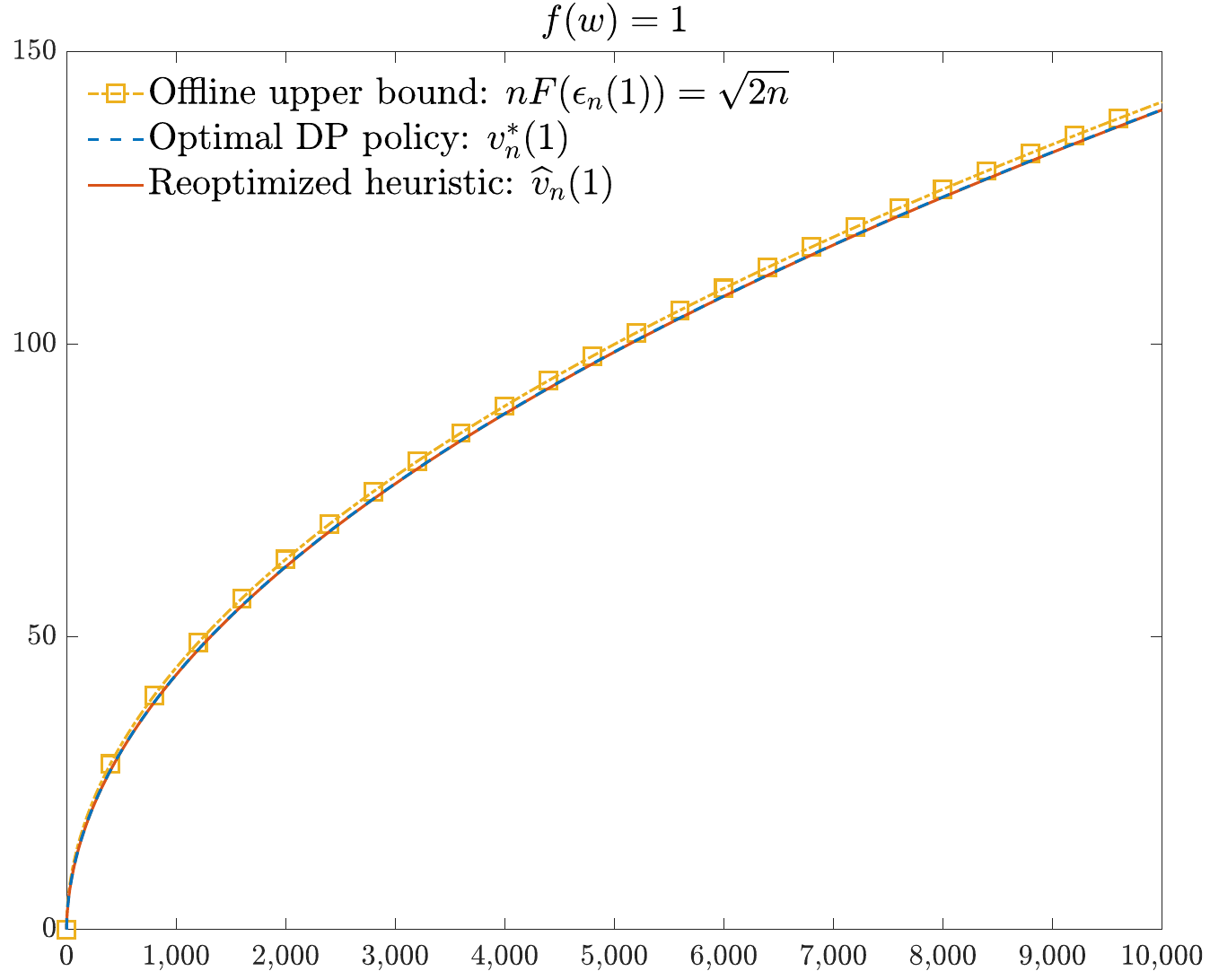}	
		\hfill
		\includegraphics[width=0.485\textwidth]{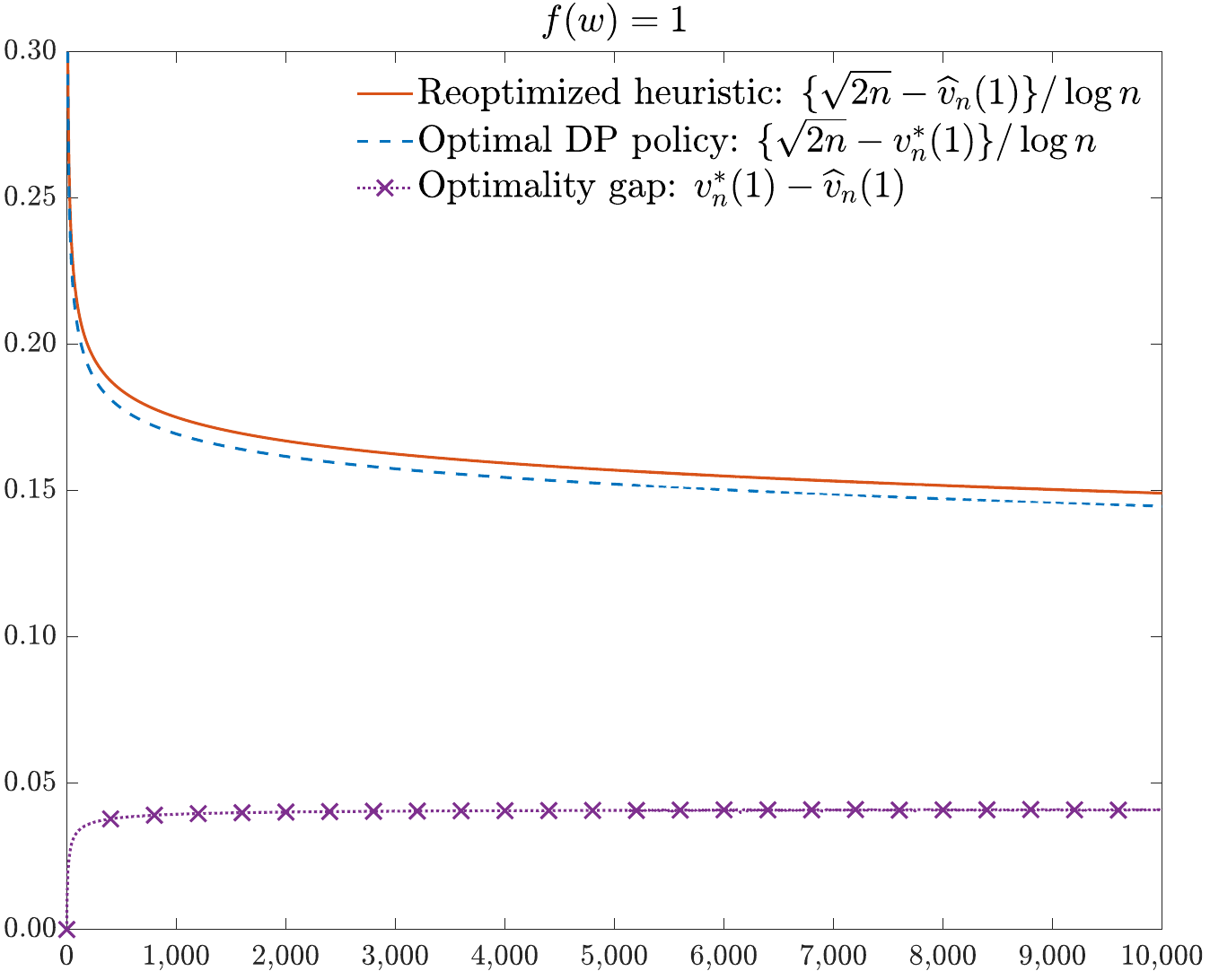}
		\\[3.5pt]
		\includegraphics[width=0.485\textwidth]{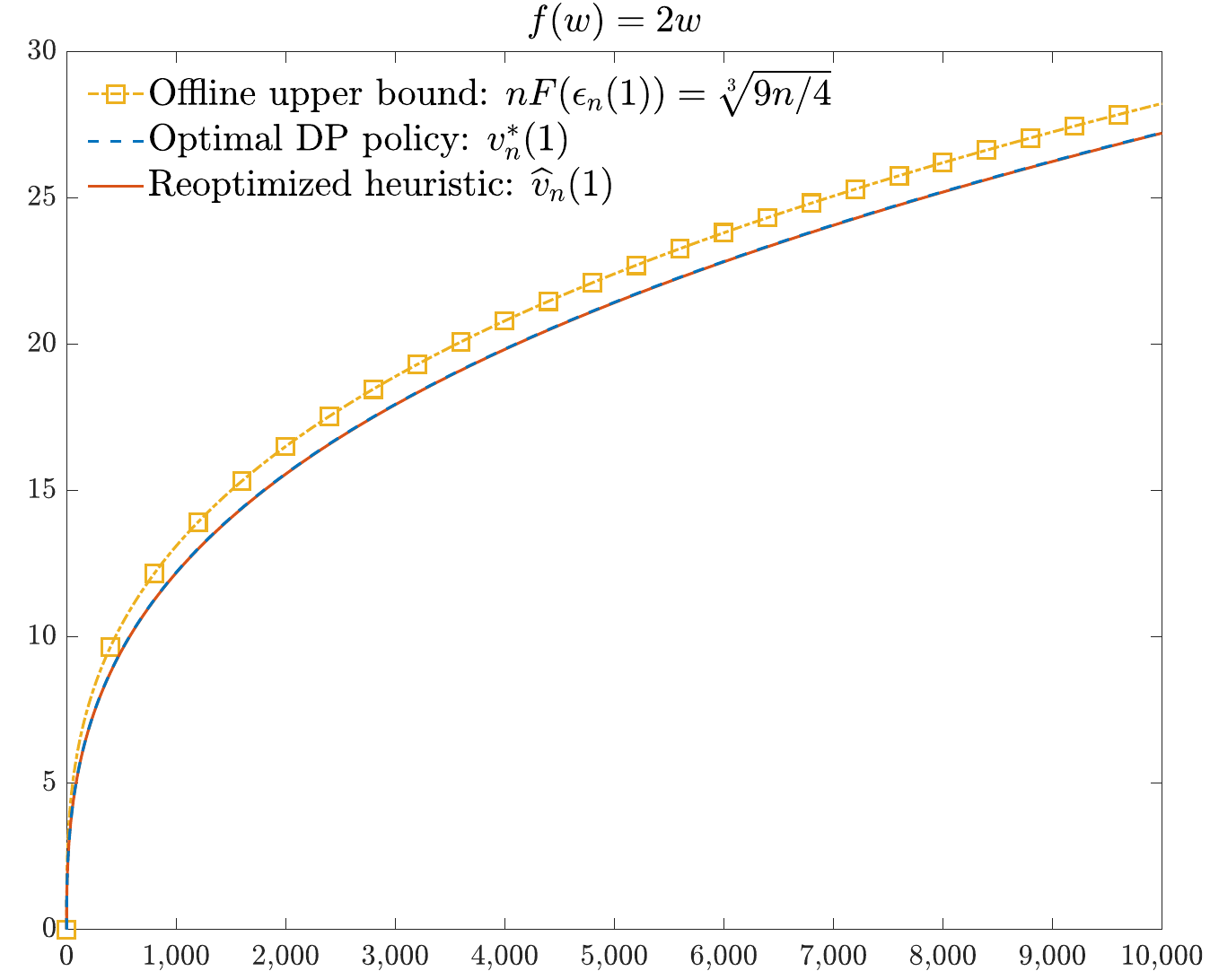}	
		\hfill
		\includegraphics[width=0.485\textwidth]{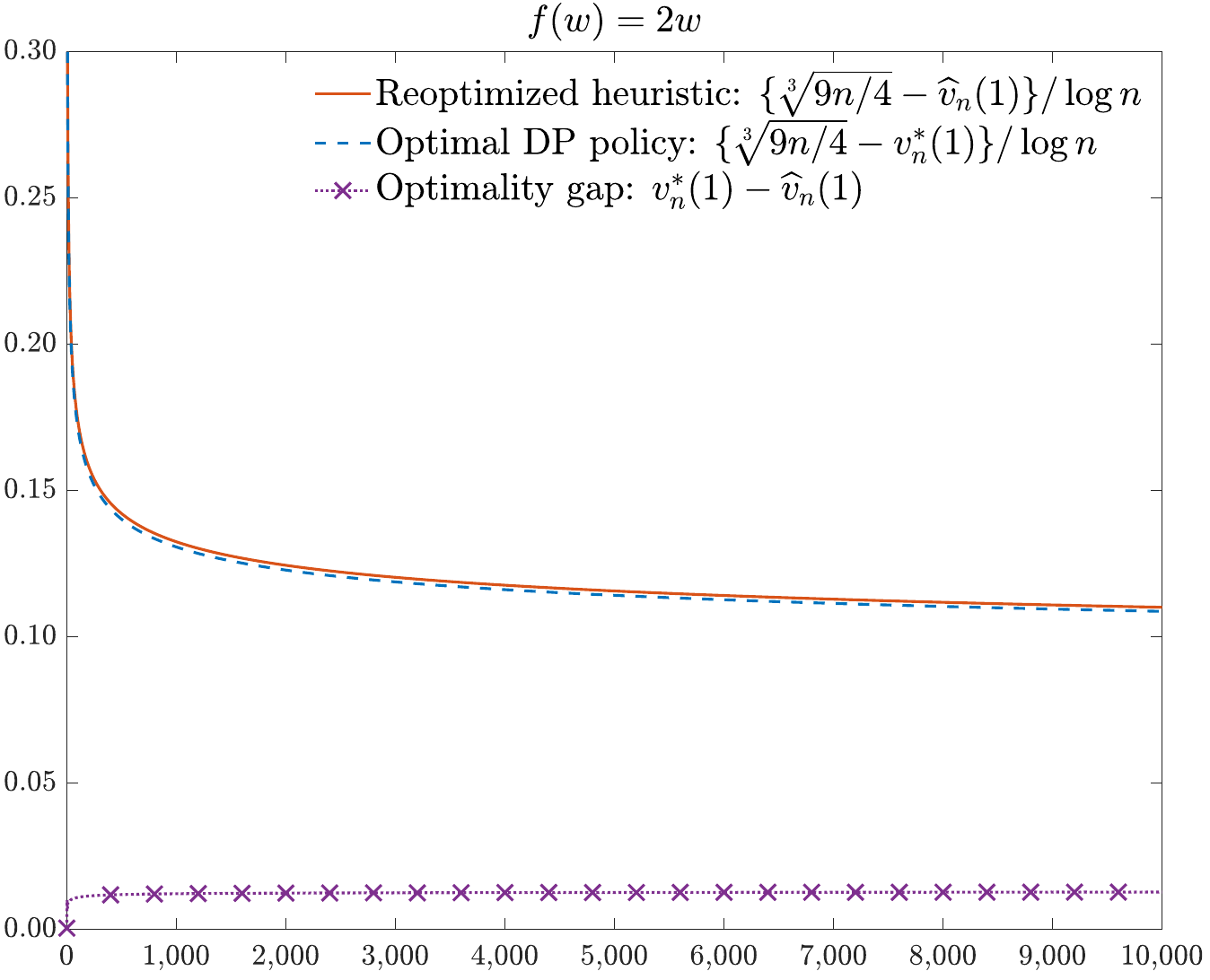}
		\\[3.5pt]
		\includegraphics[width=0.485\textwidth]{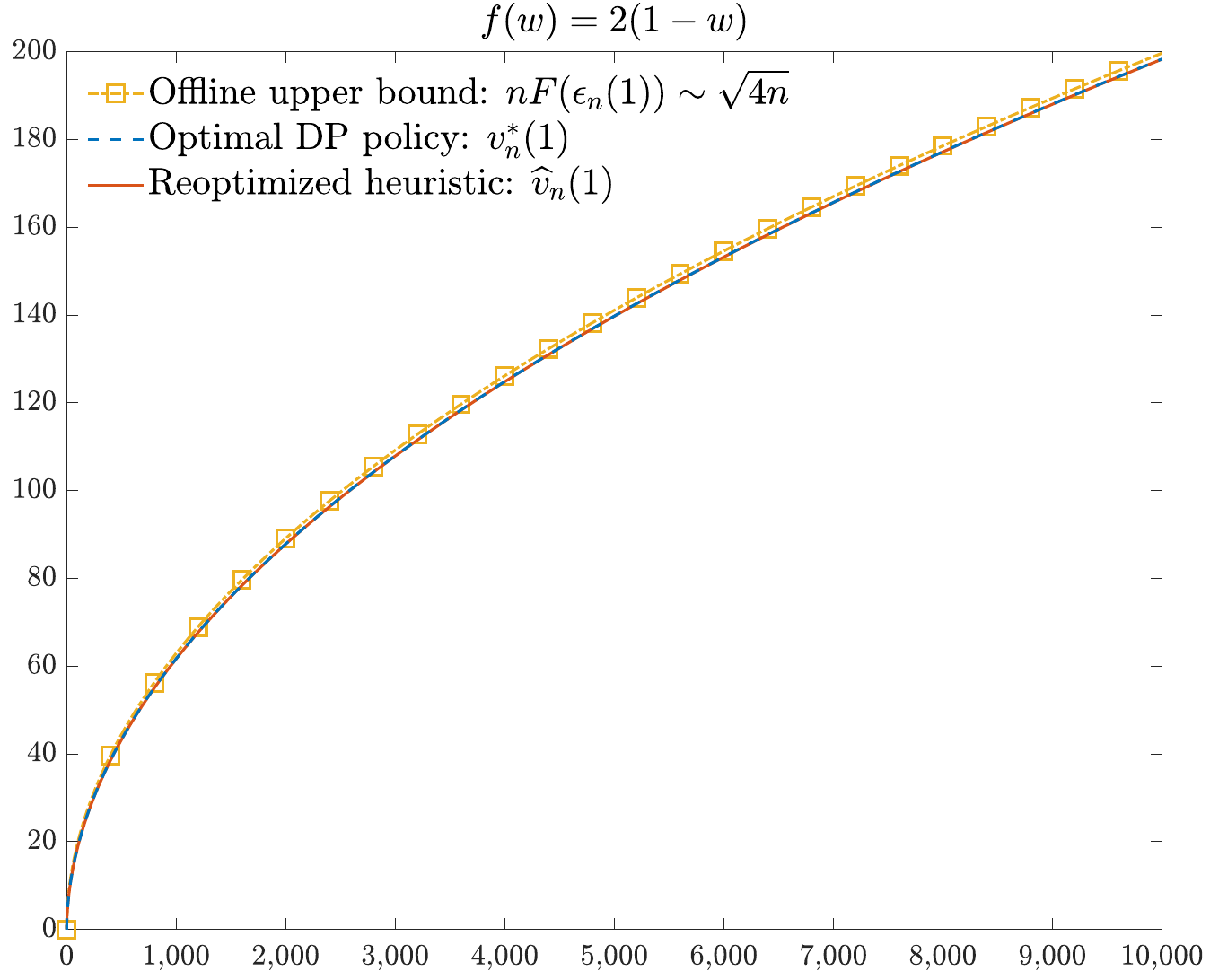}	
		\hfill
		\includegraphics[width=0.485\textwidth]{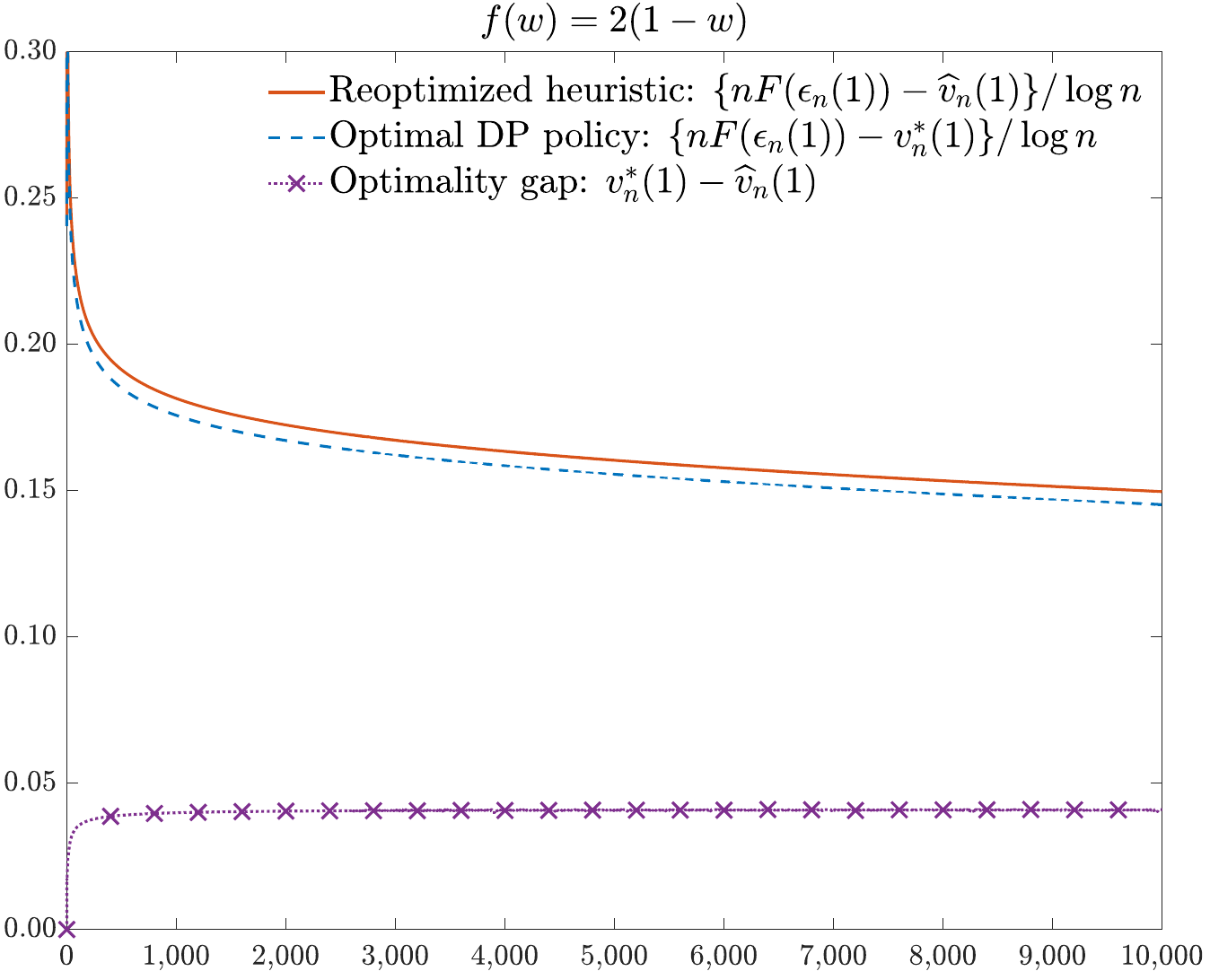}
        \\[3.5pt]
	\end{center}

	\scriptsize{\emph{Notes.}	
        The left plots display the prophet upper bound and the value functions
        of the optimal dynamic programming (DP) policy and of the reoptimized heuristic.
        The right plots show the regret bounds of the optimal policy and of the heuristic
        scaled by the logarithm of $n$, as well as the optimality gap.
        While the scaled regret bounds are bounded away from zero for large $n$,
        the optimality gap does not grow with $n$.
        Weights have densities on $(0,1)$ respectively given by $f(w)=1$, $f(w)=2w$, and $f(w)=2(1-w)$.
        Capacity $c=1$, arrival probability $p=1$, and reward $r=1$.
        Discretized state space with grid size $10^{-5}$.
        }
\end{figure}

Theorem \ref{thm:main} tells us that the regret of a dynamic and stochastic
knapsack problem is at most logarithmic in $n$, provided that the weight distribution
belongs to the typical class of Definition \ref{def:typical-class}.
While the actual order of the regret may---in principle---be smaller than
what our bound predicts, we find numerically that this is not the case.
In fact, we conjecture that the actual regret is $O(\log n)$ as $n \rightarrow \infty$
for most continuous weight distributions.

As discussed in Section \ref{se:introduction},
the work of \citet{Seksenbayev:WP2018} and \citet{GnedinSeksenbayev:WP2019}
tells us that when the capacity, the reward, and the arrival probability are all equal to one,
and the weight distribution is uniform on the unit interval, then the regret is asymptotic to $(\log n)/12$.
In this section, we numerically investigate the actual order of the regret for two other weight distributions,
while keeping the uniform as reference.

For our numerical examples, we solve the recursion \eqref{eq:v-hat}
on a discretized state space  with grid size $10^{-5}$ and obtain estimates for the reoptimized value function $\widehat{v}_n(\cdot)$
for $n \in \{1,\ldots,10000\}$ and for different distributions $F$.
On the same discretized state space and for the same weight distributions, we also solve numerically
the Bellman recursion
\begin{eqnarray}\label{eq:v-star}
	v^*_n(x)
	&=& p(1-F(x)) v^*_{n-1}(x) + p\int_0^x \max \{r+v^*_{n-1}(x-w), v^*_{n-1}(x)\} f(w) \, dw + (1-p)v^*_{n-1}(x)	\nonumber\\		
	&=& (1-pF(x)) v^*_{n-1}(x) + p\int_0^x \max \{r+v^*_{n-1}(x-w), v^*_{n-1}(x)\} f(w) \, dw
\end{eqnarray}
with the initial condition $v^*_0(x)=0$ for all $x \in [0,c]$,
and we obtain estimates for the optimal value functions $v^*_n(\cdot)$
for $n\in\{1,\ldots,10000\}$.
Finally, we simulate the average of the offline solution $R^*_n(c,p,r)$
and compare all of our numerical estimates with the prophet upper bound $nprF(\epsilon_{n p}(c))$.
Based on our numerical experiments, we observe that:
\begin{enumerate}[label=(\roman*)]
	\item The gap $nprF(\epsilon_{n p}(c)) - \E{R^*_n(c,p,r)}$ between the prophet upper bound
		  and the offline solution is bounded by a constant that does not depend on $n$
		  (see Figure \ref{fig:simulations});
	\item The regret bound $nprF(\epsilon_{n p}(c))-\widehat{v}_n(c)$ for the reoptimized
		  heuristic and the regret bound $nprF(\epsilon_{n p}(c))-v^*_n(c)$ for the optimal
		  online policy grow logarithmically with $n$ (Figure \ref{fig:value-function-and-scaled-regret}); and
	\item The optimality gap $v^*_n(c) - \widehat{v}_n(c)$ is bounded by constant
		  that is independent of $n$ (Figure \ref{fig:value-function-and-scaled-regret}).			
\end{enumerate}

In turn, our numerical experiments suggests that the regrets (rather than the regret bounds)
$\E{R^*_n(c,p,r)} - \widehat{v}_n(c)$ and $\E{R^*_n(c,p,r)} - v^*_n(c)$
respectively of the reoptimized heuristic and of the optimal online policy
are also logarithmic in $n$.
In our numerical work, we consider instances of the dynamic and stochastic knapsack problem
with reward $r=1$, arrival probability $p=1$, and capacity $c=1$.
We vary item weights by considering the three densities supported on the unit interval
given by
$f(w)=1$, $f(w)=2w$ and $f(w)=2(1-w)$ for $w \in (0,1)$.
The top left chart of Figure \ref{fig:value-function-and-scaled-regret}
plots the prophet upper bound
$nF(\epsilon_n(1)) = \sqrt{2n}$ as well as the value function of the optimal policy,
$v^*_n(1)$, and of the reoptimized heuristic, $\widehat{v}_n(1)$,
when the weight distribution is uniform on $(0,1)$.
Instead, the top right chart depicts the respective regret bounds scaled by the
logarithm of $n$, as well as the optimality gap.
In the chart we see that the scaled regret bounds (top two lines) are bounded away
from zero for large $n$, implying that the regret bounds grow logarithmically.
In contrast, the optimality gap (bottom line) appears not to grow with $n$.

The plots in the middle row of Figure \ref{fig:value-function-and-scaled-regret}
point to the same set of observations
when the weights have density $f(w)=2w \1{w \in (0,1) }$
and the prophet upper bound is $nF(\epsilon_n(1)) = \sqrt[3]{9n/4}$.
Finally, the bottom two charts of Figure \ref{fig:value-function-and-scaled-regret}
consider item weights that have density $f(w)=2(1-w)\1{w\in (0,1) }$.
In this case, the prophet upper bound
cannot be expressed in closed form, but one can show that
$nF(\epsilon_n(1)) \sim \sqrt{4n}$ as $n \rightarrow \infty$.
Nevertheless, also for this weight distribution the numerical analysis
suggests that the regrets of the optimal policy and of the heuristic
are both logarithmic in $n$,
and that the optimality gap can be bounded by a constant independent of $n$.

\section{On weight distributions with multiple types}\label{se:weights-w-multiple-types}

In this section, we discuss how our logarithm regret bound generalizes to dynamic and stochastic
knapsack problems with equal rewards and with independent random weights
that belong to one of $J < \infty$ different types.
We consider a multinomial arrival process with parameters $\boldsymbol{p} \equiv (p_0, p_1, \ldots, p_J)$
where $p_j \in (0,1]$ for all $j \in [J]$ and $p_0 = 1 - \sum_{i \in [J]} p_j  \in [0,1]$.
Here, the parameter $p_0$ represents the probability of no item arriving (or, equivalently, the arrival probability
of an item with arbitrarily large weight) and $p_j$, $j \in [J]$, is the arrival probability of an item
with weight distribution $F_j$.

Upon arrival of an item the decision maker may see the type of the item or not.
If the item types are \emph{not} released, then she only sees the arriving weights that (conditional on an arrival occurring)
are drawn from the mixture distribution
$$
	\widetilde F(w) = \frac{1}{1-p_0} \sum_{ j \in [J] }  p_j F_j(w) \quad \mbox{ for all } w \in [0,\infty).
$$
If the weight distributions $F_1, F_2, \ldots, F_J$ are all typical  (see Definition \ref{def:typical-class}),
then the mixture distribution $\widetilde F$ is also typical (see Section \ref{se:typical-class}),
and Theorem \ref{thm:main} immediately applies.

In contrast, if item types are revealed upon arrival, then the decision maker could use
the type information to make better decisions.
As we will see shortly, because the rewards are all equal, knowing the weight type of the arriving item
makes no difference. The offline solution is still given by an algorithm that sorts items according
to their realized weights (regardless of their types),
and the optimal dynamic programming policy is a threshold policy
that ignores weight types.

For the optimal offline  solution, we  can reinterpret this formulation so that items arrive according
to a Bernoulli process with arrival probability $1-p_0 = \sum_{j\in[J]} p_j$,
have rewards equal to $r$ and independent weights with distribution
given by $\widetilde F$.
The optimal offline solution $R^*_n(c, 1-p_0, r)$ is then given by the sorting algorithm \eqref{eq:offline-sort-algorithm},
so if
\begin{equation}\label{eq:epsilon-k-multi-type-mixture}
\epsilon_{k(1-p_0)}(x) = \sup \left\{ \epsilon \in[0,\infty): \int_0^{\epsilon} w \, d\widetilde F(w) \leq \frac{x}{k(1-p_0)} \right\},
\end{equation}
then Proposition \ref{prop:prophet-upper-bound} gives us that
\begin{equation}\label{eq:multi-type-upper-bound}
\E{ R^*_n(c, 1-p_0, r) } \leq n (1-p_0) r \widetilde F(\epsilon_{n (1-p_0)}(c)),
\end{equation}
and the prophet upper bound for weight distribution with multiple types follows.

To establish the independence on weight types of the optimal online solution when the rewards are all equal,
we now examine the associated Bellman equation.
We suppose that, with $k$ periods to the end of the horizon, the remaining capacity is $x \in [0,c]$,
the arriving item has weight type $j \in \{0,1,\ldots,J\}$ (with $j=0$ denoting a no arrival or, equivalently, an arrival
with arbitrarily large weight), and we let $V_k(x,j)$ be the optimal expected reward to-go given
the current state.
The optimality principle of dynamic programming then tells us that the value function $V_k(x,j)$
satisfies the Bellman recursion
\begin{equation}\label{eq:Bellman-equation-multiple-type}
	V_k(x,j) \!=\!  (1-F_j(x)) \! \sum_{\iota=0}^J p_{\iota} V_{k-1}(x,\iota) \\
                +\! \int_0^x \! \max \! \left\{ r+ \! \sum_{\iota=0}^J p_{\iota} V_{k-1}(x-w,\iota), \! \sum_{\iota=0}^J p_{\iota} V_{k-1}(x,\iota) \! \right\} dF_j(w),
\end{equation}
with the initial condition $V_0(x,j) = 0$ for all $x \in [0,c]$ and all $j \in \{0,1,\ldots,J\}$.
Here, the first summand holds because with probability $1-F_j(x)$ the arriving type-$j$
item has weight that exceeds the current knapsack capacity and the decision maker must reject it.
Thus, her expected reward to-go over the remaining $k-1$ periods is just given by the average
over types of the value functions $V_{k-1}(x,\iota)$ for $\iota \in \{0,1,\ldots, J\}$.
Instead, with probability $F_j(x)$ the arriving type-$j$ item can be selected and the decision maker chooses
the action that yields the largest expected reward to-go.
If the item has weight $w$ then its selection yields $r + \sum_{\iota=0}^J p_{\iota} V_{k-1}(x-w,\iota)$,
while its rejection gives $\sum_{\iota=0}^J p_{\iota} V_{k-1}(x,\iota)$.
By integrating this against $F_j(\cdot)$  for $w \in [0,x]$,
we obtain the second summand of \eqref{eq:Bellman-equation-multiple-type}.
The value functions $V_k(x,j)$ are monotone increasing in $x$ for each $j$ and $k$,
and one has that
$$
	H_k^*(x,j) = \sup \left\{w \in[0,x]: r +\sum_{\iota=0}^J p_{\iota} V_{k-1}(x-w,\iota) \geq \sum_{\iota=0}^J p_{\iota} V_{k-1}(x,\iota) \right\},
$$
is the optimal threshold that identifies the largest type-$j$ weight that can be selected when
the current capacity is $x$ and there are $k$ periods remaining.
Interestingly, one immediately has that $H_k^*(x,j) = H_k^*(x,\iota)$
for all $j,\iota \in \{0,1,\ldots,J\}$ since all items have the same reward $r$
and the expected rewards to-go of both actions are type independent.
Because the optimal threshold policy ignores types, we can
construct  a heuristic that has the same property and use
our earlier analysis to assess its performance.
We recall the quantity $\epsilon_{k(1-p_0)}(x)$ in \eqref{eq:epsilon-k-multi-type-mixture}
and consider the type-independent threshold
$$
	\widehat H_{k}(x,j) = \min \{x, \epsilon_{k(1-p_0)}(x)\} \qquad \text{for all } j \in \{0,1,\ldots, J\} \text{ and } x \in [0,c].
$$
If $\widehat{\pi}$ is the heuristic that uses the thresholds $\widehat{H}_n, \widehat{H}_{n-1}, \ldots, \widehat{H}_1$,
and $R^{\widehat \pi}_n(c, 1-p_0, r)$ is the total reward that $\widehat \pi$ collects,
then Proposition \ref{pr:max-residual-bound} tells us that there is a constant $1<M<\infty$ depending only on $\widetilde F$,
the arrival probability $1-p_0$, and the reward $r$ such that
\begin{equation}\label{eq:multi-type-lower-bound}
n (1-p_0) r \widetilde F(\epsilon_{n (1-p_0)}(c)) - M \log n \leq \E{ R^{\widehat \pi}_n(c, 1-p_0, r) }.
\end{equation}

If we combine the two bounds \eqref{eq:multi-type-upper-bound} and \eqref{eq:multi-type-lower-bound},
we then have the corollary below.

\begin{corollary}[Regret bound for weight distributions with multiple types]
Consider a knapsack problem with capacity $0\leq c < \infty$ and with items that arrive over $1 \leq n < \infty$ periods
according to a multinomial process with parameters $\boldsymbol{p} \equiv ( p_0, p_1, \ldots p_J)$ such that $1-p_0 = \sum_{j \in [J]} p_j$,
and where $p_0$ is the probability of no arrival.
If the items have rewards all equal to $r$ and type-dependent weights
with continuous distributions $F_1, F_2, \ldots, F_J$
and mixture (conditional on an arrival occurring) given by
$$
\widetilde F(w) = \frac{1}{1-p_0} \sum_{ j \in [J] }  p_j F_j(w) \quad \mbox{ for all } w \in [0,\infty),
$$
then
$$
 \E{R^*_n(c, 1-p_0, r) } \leq n (1-p_0) r \widetilde F(\epsilon_{n (1-p_0)}(c)).
$$
Furthermore, there is a feasible online policy $\widehat \pi$ such that
if the weights are independent and their distributions $F_1, \ldots, F_J$ belong to the typical class then there is a constant $M$
depending only on $\widetilde F$, $p_0$, and $r$ for which
$$
n (1-p_0) r \widetilde F(\epsilon_{n (1-p_0)}(c)) - M \log n \leq \E{ R^{\widehat \pi}_n(c, 1-p_0, r) }.
$$
In turn, if the weights are independent and $F_1, \ldots, F_J$ all belong to the typical class,
then we have the regret bound
$$
 \E{R^*_n(c, 1-p_0, r) }  - \E{ R^{\widehat \pi}_n(c, 1-p_0, r) }  \leq M \log n.
$$
\end{corollary}

We note here that the key assumption that makes our analysis carry over to weight
distributions with multiple types is that the rewards are all equal.
If one were to allow for type-dependent rewards, then the optimal offline solution
would \emph{not} be given by the offline-sort algorithm \eqref{eq:offline-sort-algorithm}
and the optimal online solution would \emph{not} be given by type-independent thresholds.
While one would still have a Bellman recursion analogous to \eqref{eq:Bellman-equation-multiple-type},
it is unclear how type-dependent rewards would affect our regret estimates,
and we leave this interesting open problem for future research.

\section{Conclusions and future direction}\label{se:conclusion}

In this paper we study the dynamic and stochastic knapsack problem with equal rewards
and independent random weights with common continuous distribution $F$.
We prove that---under some mild regularity conditions on the weight distribution---the regret
is, at most, logarithmic in $n$.
In particular, we show that this regret bound is attained by a reoptimized heuristic that
can be expressed explicitly and that provides a key analytical connection with the offline solution.

Two questions  stem naturally from our analysis.
The first one entails the difference in performance
between the reoptimized heuristic and the optimal online policy.
Based on our numerical experiments, we conjecture that
\begin{equation}\label{eq:conjecture-pistar-pihat}
  \max_{\pi \in \Pi(n,c,p)} \E{ R_n^{\pi}(c,r,p) }  = \E{ R_n^{\widehat \pi}(c,r,p) } + O(1)
\end{equation}
for all $n \geq 1$ and for a large class of weight distributions.
However, it is well-known that the optimal policy often lacks of desirable
structural properties, so proving \eqref{eq:conjecture-pistar-pihat}
is unlikely to be easy.
The second question has to do with the performance of the offline-sort algorithm.
Here, the numerical evidence suggests that
\begin{equation*}%\label{eq:conjecture-offline-sort}
  \E{ R^*_n(c,r,p) }  = npr \F{ \epsilon_{n p}(c) } + O(1)
\end{equation*}
for all $n \geq 1$ and most continuous weight distributions $F$.

Resolving the two conjectures above would imply that
the regret cannot be $o(\log n)$ as $n \rightarrow \infty$ for most continuous weight distributions,
and that $O(\log n)$ as $n \rightarrow \infty$ correctly quantifies the informational advantage
that the prophet has over the sequential decision maker.
This is in contrast with some other dynamic and
stochastic knapsack problems in which
the sequential decision maker does essentially as well as the prophet (see Section \ref{se:literature-review}).
It also suggests that when items have random weights,
then the design of near-optimal heuristics requires more care than usual.

\section*{Acknowledgement}

The authors are thankful to Santiago R. Balseiro, Itai Gurvich, and Yehua Wei for insightful discussions.
This material is based upon work supported by the National Science Foundation
under CAREER Award No. 1553274.

%\bibliography{biblio-knapsack}
%\bibliographystyle{abbrvnat}
%\bibliographystyle{agsm}
%\bibliographystyle{apalike}

\end{document}